\documentclass[a4paper, authoryear]{article}
\usepackage{techRepsPDF}
\usepackage[parfill]{parskip}
\usepackage{natbib}
\usepackage[ruled,vlined,linesnumbered]{algorithm2e} \SetCommentSty{textsf}
\usepackage{setspace} 
\usepackage{graphicx} 
\usepackage{placeins}
\usepackage{newproof} 
\usepackage{setspace} 
\usepackage{amsmath, amssymb}
\usepackage{xcolor, colortbl, color} \usepackage{threeparttable, longtable}
\usepackage[margin=1cm]{caption}
\usepackage[pdftitle={General Arc-flow Formulation with Graph Compression}, pdfauthor={Filipe Brandao}]{hyperref}
\usepackage{tikz} \usepackage{subfig} \usetikzlibrary{shapes.geometric}
\usetikzlibrary{shapes.multipart}
\usetikzlibrary{decorations,decorations.markings}
\usetikzlibrary{shapes,arrows,fit,calc,positioning} 
\usepackage{fullpage}
\usepackage{comment}
\usepackage{multirow}
\usepackage{amsmath}

\usepackage[utf8]{inputenc}
\usepackage[english]{babel}

\newcolumntype{$}{>{\global\let\currentrowstyle\relax}} 
\newcolumntype{^}{>{\currentrowstyle}}

\newtheorem{property}{Property}
\newtheorem{definition}{Definition} 
\newtheorem{example}{Example}

\def\HD#1#2{\vrule height #1pt depth #2pt width 0pt\relax} 
\def\up{\HD{10}{0}}
\def\down{\HD{0}{5}}

\newcommand*\centermathcell[1]{\omit\hfil$\displaystyle#1$\hfil\ignorespaces}

\newcommand{\sip}{\mathop{\mathrm{bb}}} 
\newcommand{\slp}{\mathop{\mathrm{lp}}} 
\newcommand{\sgg}{\mathop{\mathrm{gg}}} 

\newcommand{\zlp}{z^{\slp}} 
\newcommand{\zip}{z^{\sip}}

\newcommand{\tlp}{t^{\slp}}
\newcommand{\tgg}{t^{\sgg}}

\newcommand{\sbb}{\mathop{\mathrm{bb}}} 
\newcommand{\nbb}{n^{\sbb}}
\newcommand{\spp}{\mathop{\mathrm{pp}}} 
\newcommand{\tpp}{t^{\spp}}
\newcommand{\tip}{t^{\sip}}
\newcommand{\stot}{\mathop{\mathrm{tot}}} 
\newcommand{\ttot}{t^{\stot}}

\newcommand{\BigO}{\mathcal{O}}
\newcommand{\vS}{\textsc{s}} 
\newcommand{\vT}{\textsc{t}}

\newcommand{\Deg}{\mbox{deg}} 
\newcommand{\Adj}{\mbox{adj}}

\newcommand{\XVEC}{\mathbf{x}}

\def\svdots{\lower1.5pt\hbox{$\smash\vdots$}}
\def\sddots{\lower1.5pt\hbox{$\smash\ddots$}}

\def\ptitle{
Bin Packing and Related Problems:\\
General Arc-flow Formulation with\\ Graph Compression
}
\def\pnumber{DCC-2013-08}
\def\pauthor{
Filipe Brand\~ao\\
{\small INESC TEC and Faculdade de Ci\^{e}ncias,
Universidade do Porto, Portugal}\\
{\small\texttt{fdabrandao@dcc.fc.up.pt}}
\ \\ 
\ \\
Jo\~ao Pedro Pedroso\\
{\small INESC TEC and Faculdade de Ci\^{e}ncias,
Universidade do Porto, Portugal}\\
{\small\texttt{jpp@fc.up.pt}}
}

\title{\textbf{\ptitle}}
\author{\pauthor}
\date{26 September, 2013}

\graphicspath{{charts/}}
\graphicspath{{figures_flow/}}

\begin{document}

\trtitle{\ptitle}
\trauthor{\pauthor}
\trnumber{\pnumber}

\mkcoverpage

\maketitle

\begin{abstract}
We present an exact method, based on an arc-flow formulation
with side constraints, for solving bin packing and cutting stock problems ---
including multi-constraint variants --- by simply representing
all the patterns in a very compact graph.
Our method includes a graph compression algorithm
that usually reduces the size of the underlying graph 
substantially without weakening the model.
As opposed to our method, which provides strong models,
conventional models are usually highly symmetric and provide very weak
lower bounds.

Our formulation is equivalent to Gilmore and Gomory's, thus providing
a very strong linear relaxation.  However, instead of using
column-generation in an iterative process, the method constructs a
graph, where paths from the source to the target node represent
every valid packing pattern. 

The same method, without any problem-specific parameterization,
was used to solve a large variety of instances from several different cutting and packing
problems.
In this paper, we deal with
vector packing, graph coloring,
bin packing, cutting stock, cardinality constrained bin packing,
cutting stock with cutting knife limitation, cutting stock 
with binary patterns, bin packing with conflicts, 
and cutting stock with binary patterns and forbidden pairs.
We report computational results obtained with many benchmark test data sets, 
all of them showing a large advantage of this formulation with respect 
to the traditional ones.
\end{abstract}

\section{Introduction}

The bin packing problem (BPP) is a combinatorial NP-hard problem (see, e.g.,
~\citealt{Garey:1979:CIG:578533}) in which objects of different
volumes must be packed into a finite number of bins, each with
capacity $W$, in a way that minimizes the number of bins used.
In fact, the BPP is strongly NP-hard (see, e.g., \citealt{Garey:1978:SNR:322077.322090}) 
since it remains so even when all of its numerical parameters are bounded by a polynomial in the length of the input.
Therefore, the BPP cannot even be solved in pseudo-polynomial time unless P = NP.
Besides being strongly NP-hard, the BPP is also hard to approximate within $3/2 -
\varepsilon$. If such approximation exists, one could partition $n$~non-negative 
numbers into two sets with the same sum in polynomial
time.  This problem --- called the number partition problem --- 
could be reduced to a bin packing problem
with bins of capacity equal to half of the sum of all the numbers.
Any approximation better than $3/2-\varepsilon$ of the optimal value could be used to find a
perfect partition, corresponding to a packing in $\lfloor 2(3/2 - \varepsilon)\rfloor = 2$ bins.
However, the number partition problem is known to be NP-hard.
Concerning the heuristic solution of the BPP,
\cite{NAV:NAV3220410409} showed that the first-fit decreasing and the best-fit
decreasing heuristics have an absolute performance ratio of~$3/2$,
which is the best possible absolute performance ratio for the bin
packing problem unless P=NP.

The BPP can be seen as a special case of the cutting
stock problem (CSP). In this problem there is a number of rolls of paper of
fixed width waiting to be cut for satisfying demand of different
customers who want pieces of various widths. 
Rolls must be cut in such a way that waste is
minimized.  Note that, in the paper industry, solving this problem to
optimality can be economically significant; a small improvement in
reducing waste can have a huge impact in yearly savings.

There are many similarities between BPP and CSP.  
However, in the CSP, 
the items of equal size (which are usually ordered in large
quantities) are grouped into orders with a required level of demand,
while in the BPP the demand for a given size is
usually close to one. 
According to \cite{Waescher}, cutting stock problems                
are characterized by a weakly heterogeneous assortment of small items,
in contrast with bin packing problems,
which are characterized by a strongly heterogeneous
assortment of small items.

The $p$-dimensional vector bin packing problem ($p$D-VBP),
also called general assignment problem by some authors,
is a generalization of bin packing with multiple constraints.
In this problem, we are required to pack $n$~items
of $m$~different types, represented  by $p$-dimensional vectors, 
into as few bins as possible.
In practice, this problem models, for example,
static resource allocation problems where
the minimum number of servers with known capacities 
is used to satisfy a set of services with known demands.

The method presented in this paper allows solving several cutting and packing
problems through reductions to vector packing.
The reductions are made by defining a matrix of weights,
a vector of capacities and a vector of demand.
Our method builds very strong integer programming models
that can usually be easily solved using any state-of-the-art mixed integer programming
solver.
Computational results obtained with many benchmark test data sets 
show a large advantage of our method with respect 
to traditional ones.

The remainder of this paper is organized as follows. 
Section~\ref{sec:exact_previous_work} 
gives account of previous approaches with exact methods to bin packing and
related problems. 
Our method is presented in Section~\ref{sec:new method}.
Computational results are
presented in Section~\ref{sec:results}, and 
Section~\ref{sec:conclusions} presents the conclusions.

\section{Previous work on exact methods}
\label{sec:exact_previous_work}

In this section, we will give account of previous approaches
with exact methods to bin packing and related problems. 
We will introduce the assignment-based Kantorovich's formulation
and the pattern-based Gilmore and Gomory's formulation.
Both formulations were initially proposed for the standard CSP,
but they can be easily generalized for multi-constraint variants.
\cite{Valerio:02} provides an excellent survey on integer programming
models for bin packing and cutting stock problems.  Here we will just
look at the most common and straightforward approaches.
 
We will also introduce Val\'erio de Carvalho's arc-flow formulation,
which is equivalent to Gilmore and Gomory's formulation.
Gilmore and Gomory's model provides a very strong linear
relaxation, but it is potentially exponential in the number
of variables with respect to the input size; 
even though Val\'erio de Carvalho's model is also potentially exponential, 
it is usually much smaller, being pseudo-polynomial in terms of decision
variables and constraints.
In both models, we consider every valid packing pattern.  
However, in Val\'erio de Carvalho's model, patterns are 
derived from paths in a graph, whereby the model is
usually much smaller.

\subsection{Kantorovich's formulation}

\cite{Kantorovich:1} introduced the following mathematical
programming formulation for the CSP:
\begin{alignat}{3}
  & \mbox{minimize }   && \sum_{k=1}^{K} y_{k}  \label{eq:kantorovich1}\\  
  & \mbox{subject to } \qquad&& \sum_{k = 1}^{K} x_{ik} \geq b_i,  && i=1,\ldots,m, \label{eq:kantorovich2}\\
  &         && \sum_{i = 1}^{m} w_{i} x_{ik} \leq W y_k, \qquad    && k=1,\ldots,K,  \label{eq:kantorovich3}\\
  &         && y_{k} \in \{0,1\},                       && k=1,\ldots,K,  \label{eq:kantorovich4}\\
  &         && x_{ik} \geq 0, \mbox{ integer},          && i=1,\ldots,m, \; k=1,\ldots,K,  \label{eq:kantorovich5}
\end{alignat}
where $K$ is a known upper bound to the number of rolls needed, $m$~is
the number of different item sizes, $w_i$~and~$b_i$ are the weight and
demand of item~$i$, and $W$~is the roll length.  The variables are~$y_k$, 
which is 1 if roll~$k$ is used and 0 otherwise, and $x_{ik}$,
the number of times item~$i$ is cut in the roll~$k$.

This model can be generalized for multi-constraint variants using
multiple knapsack constraints instead of just one. 
In the $p$-dimensional  case, the lower bound provided by the
linear relaxation of this kind of assignment-based models 
is~$1/(p+1)$ (see, e.g., \citealt{Caprara:1998:PIF:303366.303371}).
Even for the one-dimensional case, the lower bound provided by this model approaches
$1/2$~of the optimal solution in the worst case.
This is a drawback of this model, as good quality lower bounds are
vital in branch-and-bound procedures.  Another drawback is due to the
symmetry of the problem, which makes this model very inefficient in
practice.

Dantzig-Wolfe decomposition is an algorithm for solving linear programming
problems with a special structure (see, e.g., \citealt{Dantzig_Wolfe_1960}). 
It is a powerful tool that can be used to obtain
models for integer and combinatorial optimization problems with stronger linear
relaxations.
\cite{Vance:1998:BAO:290008.290010} applied a Dantzig-Wolfe
decomposition to model
(\ref{eq:kantorovich1})-(\ref{eq:kantorovich5}), keeping constraints
(\ref{eq:kantorovich1}), (\ref{eq:kantorovich2}) in the master
problem, and the subproblem being defined by the integer solutions to
the knapsack constraints (\ref{eq:kantorovich3}).  
Vance showed that when all the rolls have the same width, the reformulated model is
equivalent to the classical Gilmore-Gomory's model.

\subsection{Gilmore-Gomory's formulation}

\cite{gomory1} proposed the following model for the CSP.
A combination of orders in the width of the
roll is called a cutting pattern.  Let column vectors $a^j = (a_1^j,
\ldots, a_m^j)^{\top}$ represent all possible cutting patterns~$j$. The
element~$a_i^j$ represents the number of items of width $w_i$ obtained
in cutting pattern~$j$.  Let~$x_j$ be a decision variable that
designates the number of rolls to be cut according to cutting pattern~$j$.  
The CSP can be modeled in terms of these
variables as follows:
\begin{alignat}{3}
  & \mbox{minimize }          && \sum_{j \in J} x_j \label{eq:cuttingstock1}\\  
  & \mbox{subject to } \qquad && \sum_{j \in J} a_i^j x_j \geq b_i, \qquad && i=1,\ldots,m, \label{eq:cuttingstock2}\\
  &                           && x_{j} \geq 0, \mbox{ integer},    \qquad && \forall j \in J, \label{eq:cuttingstock3}
\end{alignat}
where $J$~is the set of valid cutting patterns that satisfy:
\begin{equation}
\sum_{i = 1}^{m} a_i^j  w_i \leq W \mbox{ and }  a_i^j \geq 0, \mbox{ integer}. \label{eq:cuttingstock4}
\end{equation}
Since constraints~(\ref{eq:cuttingstock4})
just accept integer linear combinations of items, 
the search space of the continuous relaxation is reduced 
and the lower bound provided by its linear relaxation is stronger
when compared with Kantorovich's formulation.

Since it may be impractical to enumerate all the columns
in the previous formulation, \cite{gomory2} proposed column
generation.
Let Model~(\ref{eq:cuttingstock1})-(\ref{eq:cuttingstock3})
be the restricted master problem.
At each iteration of the column-generation process, a subproblem
is solved and a column (pattern) is introduced in the restricted master problem
if its reduced cost is strictly less than zero.
The subproblem, which is a knapsack problem, is the following:
\begin{alignat}{3}
  & \mbox{minimize}          && 1-\sum_{i = 1}^{m} c_i a_i \label{eq:subprob1}\\  
  & \mbox{subject to } \qquad && \sum_{i = 1}^{m} w_i a_i \leq W\label{eq:subprob2}\\
  &                           && a_i \geq 0, \mbox{ integer},    \qquad &&  i=1,\ldots,m,  \label{eq:subprob3}
\end{alignat}
where $c_i$~is the shadow price of the demand constraint of item~$i$
obtained from the solution of the linear relaxation of the restricted master problem,
and $a = (a_1, \ldots, a_m)$ is a cutting pattern 
whose reduced cost is given by the objective function.
This model can be easily generalized for multi-constraint variants
by using multiple knapsack constraints in the subproblem
instead of (\ref{eq:subprob2}).

\subsection{Val\'erio de Carvalho's arc-flow formulation}
\label{sec:valerio}

Among other methods for solving BPP and CSP exactly, 
one of the most important is the arc-flow formulation with side constraints 
of \cite{Valerio:01}.
This model has a set of flow conservation constraints and a set of 
demand constraints to ensure that
the demand of every item is satisfied. The corresponding path-flow formulation is
equivalent to the classical Gilmore-Gomory's formulation.

Consider a bin packing instance with bins of capacity~$W$ and 
items of sizes $w_1,w_2,\ldots,w_m$ with demands $b_1, b_2, \ldots, b_m$, respectively.
The problem of determining a valid solution
to a single bin can be modeled as the problem of finding a path in a
directed acyclic graph $G=(V,A)$ with $V = \{0,1,2,\ldots,W\}$ and $A =
\{(i,j)\ |\ j-i=w_d, \mbox{ for } 1 \leq d \leq m \mbox{ and } 0 \leq i
< j \leq W\}$, meaning that there exists an arc between two vertices
$i$ and $j>i$ if there are items of size $w_d = j - i$.  The number
of vertices and arcs are bounded by $\BigO(W)$ and $\BigO(m W)$, respectively.  
Additional arcs $(k,k+1)$, for $k = 0,\ldots,W-1$, are included for representing
unoccupied portions of the bin.

In order to reduce the symmetry of the solution space 
and the size of the model, Val\'erio de Carvalho introduced some 
rules. The idea is to consider
only a subset of arcs from~$A$.
If we search for a solution in which the items are ordered by decreasing
values of width, only paths in which items appear according
to this order must be considered.
Note that different paths that include the same set of items
are equivalent and we just need to consider one of them.

\begin{figure}[h!tbp]
\caption{Graph corresponding to Example~\ref{ex:example1}.\label{fig:graph2}}
  \centering
  \includegraphics[scale=1]{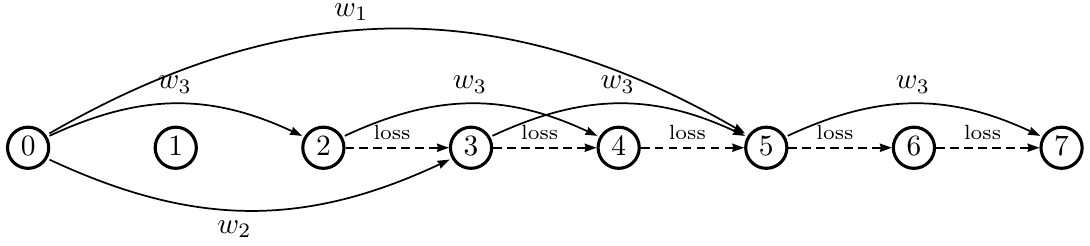}
\end{figure}

\begin{example}
\label{ex:example1}
Figure~\ref{fig:graph2} shows the graph associated with an instance
with bins of capacity~$W = 7$ and items of sizes 5, 3, 2 with demands
3, 1, 2, respectively.
\end{example}

The BPP and the CSP are thus equivalently formulated as that of
determining the minimum flow between vertex~$0$ and vertex~$W$, with
additional constraints enforcing the sum of the flows in the arcs of
each order to be greater than or equal to the corresponding demand.
Consider decision variables~$x_{ij}$ (associated with arcs~$(i,j)$ defined above) 
corresponding to the number of items of size~$j - i$ placed in any 
bin at a distance of $i$~units from the 
beginning of the bin. 
A variable~$z$, representing the number of bins required, aggregates the
flow in the graph, and can be seen as a feedback arc from vertex~$W$
to vertex~$0$.  The model is as follows:
\begin{alignat}{3}
  & \mbox{minimize }   && z  \label{eq:valerio1}\\  
  & \mbox{subject to } \qquad&& \sum_{(i,j)\in A} x_{ij} - \sum_{(j,k) \in A}x_{jk} = && \left\{ \begin{array}{rl}
                -z & \mbox{if }  j = 0, \\
                z  & \mbox{if }  j = W, \\
                0  & \mbox{for } j = 1,\ldots,W-1,\\
                \end{array}\right. \label{eq:valerio2} \\
      &         && \sum_{(k,k+w_i) \in A} x_{k,k+w_i} \geq b_i, \quad && i=1,\ldots,m, \label{eq:valerio3} \\
      &         && x_{ij} \geq 0, \;  \mbox{integer},              && \forall (i,j) \in A. \label{eq:valerio4}
\end{alignat}

\cite{Valerio:01} developed a branch-and-price
procedure that combines column-generation and
branch-and-bound.  At each iteration, the subproblem generates a set
of variables, which altogether correspond to an attractive valid packing
for a single bin.  

\section{General arc-flow formulation}
\label{sec:new method}

In this section, we propose a generalization of Val\'erio de Carvalho's arc-flow
formulation.
Val\'erio de Carvalho's graph can be seen as the dynamic programming search
space of the underlying one-dimensional knapsack problem.
The vertices of the graph can be seen as states
and, in order to model multi-constraint knapsack problems,
we just need to add more information to them.

\begin{figure}[h!tbp]
\caption{Another possible graph for a bin packing instance.\label{fig:bpgraph}}

  \centering
      \begin{minipage}[c]{0.64\linewidth}
        \includegraphics[scale=1]{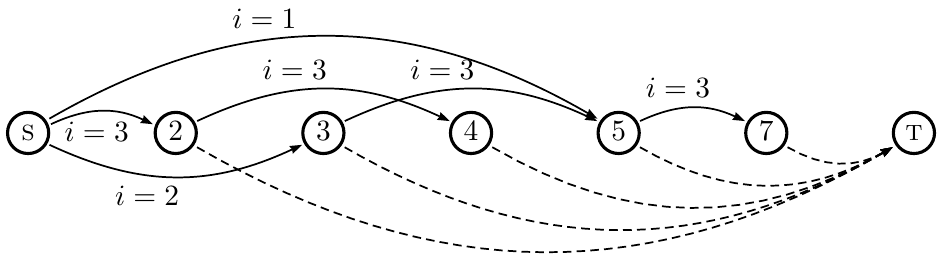}
      \end{minipage}  
      \begin{minipage}[c]{0.1\linewidth}      
        \begin{tabular}{ccc}
                $i$ & $w_i$ & $b_i$ \\
                \hline
                \up\down
                1 & 5 & 3\\
                2 & 3 & 1\\
                3 & 2 & 2\\
                \hline
                \multicolumn{3}{c}{\up\down $W = 7$}\\
        \end{tabular}\\ 
      \end{minipage} 
  
\caption*{\footnotesize 
In this graph, $V = \{\vS, 2, 3, 4, 5, 7, \vT\}$ is the set of vertices,
$A =$ \{(\vS, 5, 1), (\vS, 3, 2), (\vS, 2, 3), (2, 4, 3), (3, 5, 3), (5, 7, 3),
(2, \vT, 0), (3, \vT, 0), (4, \vT, 0), (5, \vT, 0), (7, \vT, 0)\} is the set of arcs.
}
\end{figure}

Figure~\ref{fig:bpgraph} shows an alternative graph to model 
the bin packing instance of Example~\ref{ex:example1} (Section~\ref{sec:valerio}).
Item weights are sorted in decreasing order.
In this graph, an arc~$(u, v, i)$ corresponds to an arc between nodes~$u$ and~$v$ 
associated with items of weight~$w_i$.
Note that, for each pair of nodes~$(u, v)$, multiple arcs associated with different items
are allowed.
The dashed arcs are the loss arcs that connect every node (except the source ($\vS$))
to the target ($\vT$).
Since the loss arcs connect the nodes directly to the target (instead
of connecting consecutive nodes) it may not be necessary to have a node for
every integer value less than or equal to the capacity.
Note that each path between~$\vS$ and~$\vT$ in this graph corresponds to a
valid packing pattern for the Example~\ref{ex:example1} and
all the valid patterns are represented.
In this graph, a node label~$w^{\prime}$ means that every sub-pattern from the source
to this node has at most weight~$w^{\prime}$.

\FloatBarrier

We propose a generalization of Val\'erio de Carvalho's model, 
still based on representing the packing patterns by means of flow
in a graph, though we assign a more general meaning to vertices and arcs.
The formulation is the following:
\begin{alignat}{3}
  & \mbox{minimize }   && z  \label{eq:new1}\\  
  & \mbox{subject to } \qquad && \sum_{(u,v,i)\in A:v=k} \hspace{-5mm}f_{uvi} \hspace{3mm}-\hspace{-3mm} \sum_{(v,r,i) \in A:v=k} \hspace{-5mm}f_{vri} = &&
  \left\{ 
    \begin{array}{rl}
      -z & \mbox{if } k = \vS, \\
      z & \mbox{if } k = \vT,\\
      0 & \mbox{for } k \in V \setminus \{\vS, \vT \},\\
    \end{array} \label{eq:new2}
  \right.\\    
      &&         & \sum_{(u, v, j) \in A:j=i} \hspace{-5mm}f_{uvj} \geq b_i, && i \in \{1,\ldots,m\} \setminus J, \label{eq:new3}\\ 
      &&         & \sum_{(u, v, j) \in A:j=i} \hspace{-5mm}f_{uvj} = b_i, && i \in J,  \label{eq:new4}\\
      &&         & f_{uvi} \leq b_i, && \forall (u,v,i) \in A, \mbox{ if } i \neq 0, \label{eq:new5}\\
      &&         & f_{uvi} \geq 0, \mbox{ integer}, && \forall (u,v,i) \in A, \label{eq:new6}      
\end{alignat}
where $m$~is the number of different items,
$b_i$~is the demand of items of weight $w_i$,
$V$~is the set of vertices,
$\vS$~is the source vertex and $\vT$~is the target;
$A$~is the set of arcs, where
each arc has three components~$(u, v, i)$
corresponding to an arc between nodes~$u$ and~$v$
that contributes to the demand of items of weight~$w_i$;
arcs~$(u,v,i=0)$ are the loss arcs;
$f_{uvi}$~is the amount of flow along the arc~$(u, v, i)$;
and $J \subseteq \{1,\ldots,m\}$ is a subset of items whose demands 
are required to be satisfied exactly for efficiency purposes.
For having tighter constraints, one may set $J = \{i = 1,\ldots, m\ |\ b_i = 1\}$ 
(we have done this in our experiments).

In Val\'erio de Carvalho's model, a variable~$x_{ij}$
contributes to an item with weight~$j-i$. In our model,
a variable~$f_{uvi}$ contributes to items of weight $w_i$;
the label of~$u$ and~$v$ may have no direct relation to the item's weight.
This new model is more general; Val\'erio de Carvalho's model is a sub-case, 
where an arc between nodes~$u$ and~$v$ can only contribute
to the demand of an item of weight~$v-u$.
As in Val\'erio de Carvalho's model, each arc can only contribute to an item,
but the new model has several differences with respect to the original formulation:
\begin{itemize}
\item nodes are more general (e.g., they can encompass multiple dimensions);
\item there may be more than one arc between two vertices (multigraph);
\item demands in general may be satisfied with excess but
for some items they are required to be satisfied exactly (this allows, 
for example, to take advantage of special ordered sets of type 1
when requiring the demands of items with demand one to be satisfied exactly);
\item arcs have upper bounds equal to the total demand of the associated item
(which allows excluding many feasible solutions that would exceed the demand);
\item arc lengths are not tied to the corresponding item weight (i.e., $(u,v,i) \in A$ even if $v-u \neq w_i$).
\end{itemize}

More details on algorithms for graph construction are
given in the following sections.
Using this model it is possible to use more general graphs, but we
always need to ensure that it is a directed acyclic graph whose paths
from~$\vS$ to~$\vT$ correspond to every valid packing pattern
of the original problem. One of the properties of this model
is the following.
\begin{property}[equivalence to the classical Gilmore-Gomory model]\label{equivalence_gomory2}
For a graph with all valid packing patterns represented as paths 
from $\vS$ to $\vT$, model (\ref{eq:new1})-(\ref{eq:new6})  is equivalent to the classical 
Gilmore-Gomory model~(\ref{eq:cuttingstock1})-(\ref{eq:cuttingstock4}) 
with the same patterns as the ones 
obtained from paths in the graph.
\end{property}
\begin{proof}  Extending Val\'erio de Carvalho's proof, we apply
Dantzig-Wolfe decomposition to model~(\ref{eq:new1})-(\ref{eq:new6})
keeping (\ref{eq:new1}), (\ref{eq:new3}) and (\ref{eq:new4}) in the master problem and
(\ref{eq:new2}), (\ref{eq:new5}) and (\ref{eq:new6}) in the subproblem. As the subproblem
is a flow model that will only generate patterns resulting
from paths in the graph, we can substitute (\ref{eq:new2}), (\ref{eq:new5}) and (\ref{eq:new6}) 
by the patterns and obtain the classical model.  
From this equivalence follows that lower bounds
provided by both models are the same when the same set of patterns is considered.
The equality constraints (\ref{eq:new4}) and the upper bound on 
variable values (\ref{eq:new5})
have no effect on the lower bounds, since for every optimal solution
with some excess there is a solution with the same objective
value that satisfies the demand exactly;
this solution can be obtained by replacing the use of some
patterns by other patterns that do not include the items whose
demands are being satisfied with excess (recall that
every valid packing pattern is represented in the graph).
\end{proof}

After having the solution of the arc-flow integer optimization model,
we use a flow decomposition algorithm to obtain the corresponding packing solution.
Flow decomposition properties (see, e.g., \citealt{ahuja-magnanti-orlin-93})
ensure that non-negative flows can be represented by paths and cycles.
Since we require an acyclic graph, any valid flow can be decomposed
into directed paths connecting the only excess node (node \vS)
to the only deficit node (node \vT).

Note that by requiring the demands of some items to be satisfied exactly
and by introducing upper bounds on variable values,
the model may sometimes become harder to solve. For instance,
if the demand of a very small item is required to be satisfied exactly,
many optimal solutions will probably be excluded (unless the optimal solution
has waste smaller than the item size). In some instances,
the solution can be obtained more quickly by choosing carefully the variable upper bounds
and the set of items whose demand must be satisfied exactly.
Overall, our choices regarding these two aspects proved
to work very well, as we show in Section~\ref{sec:results}.

\subsection{$p$-dimensional vector packing graphs}
\label{sec:multidimgraphs}

\begin{figure}[h!tbp]
\caption{Graph associated with Example~\ref{ex:example1},
but now with cardinality limit 3.\label{fig:cbpgraph}}
  
  \centering
      \begin{minipage}[c]{0.70\linewidth}
        \includegraphics[scale=1]{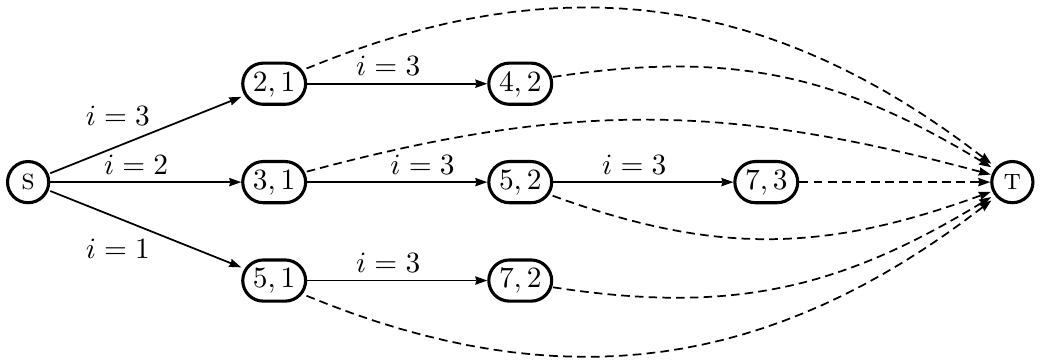}
      \end{minipage}  
      \begin{minipage}[c]{0.2\linewidth}      
        \begin{tabular}{cccc}
                $i$ & $w_i^1$ & $w_i^2$ & $b_i$ \\
                \hline 
                \up\down        
                1 & 5 & 1 & 3\\
                2 & 3 & 1 & 1\\
                3 & 2 & 1 & 2\\
                \hline								 											
                \multicolumn{4}{c}{\up\down $W^1 = 7$; $W^2 = 3$}\\
        \end{tabular}\\ 
      \end{minipage} 

\caption*{\footnotesize
In this graph, a node label $(w^{\prime} ,c^{\prime})$ means that any path from the source
will reach the node with at most $c^{\prime}$ items whose
sum of weights is at most $w^{\prime}$.
}
\end{figure}

In the $p$-dimensional vector packing problem, items
are characterized by weights in several dimensions.
For each dimension~$d$, let~$w_i^d$ be the weight of item~$i$ and $W^d$~the bin capacity.
In the one-dimensional case, arcs associated with items of weight~$w_i^1$ 
lie between vertices~$(a)$ and~$(a+w_i^1)$.
In the multi-dimensional case, arcs associated with items of weight $(w_i^1,w_i^2,\ldots,w_i^p)$ 
lie between vertices
$(a^1,a^2,\ldots,a^p)$ and $(a^1+w_i^1,a^2+w_i^2,\ldots,a^p+w_i^p)$.
Figure~\ref{fig:cbpgraph} shows the graph associated with 
Example~\ref{ex:example1}, but now with a second dimension
limiting cardinality to 3.

In the one-dimensional case, the number of vertices and arcs in the arc-flow
formulation is bounded by $\BigO(W)$ and $\BigO(mW)$, respectively,
and thus graphs are usually reasonably small.
However, in the multi-dimensional case, the number of vertices and arcs are bounded by 
$O(\prod_{d=1}^{p} (W^d+1))$ and $O(m\prod_{d=1}^{p} (W^d+1))$, respectively.
Despite the possible intractability indicated by these bounds,
the graph compression method that we present in this paper usually
leads to reasonably small graphs even for very hard instances
with hundreds of dimensions.

\begin{definition}[Order]\label{def:order}
Items are sorted in decreasing order by the sum of normalized weights
($\alpha_i = \sum_{d=1}^{p} x_i^d/W^d$), using decreasing
lexicographical order in case of a tie.
\end{definition}

\begin{figure}[h!tbp]
\caption{Graph construction example.\label{fig:graphalg1}}
\centering\fbox{\begin{minipage}[b]{0.80\textwidth}
Consider a cardinality constrained bin packing instance 
with bins of capacity~7, cardinality limit~3,
and items of sizes~5, 3, 2 with demands~3, 1, 2, respectively.
This instance corresponds to adding cardinality limit~3 to
Example~\ref{ex:example1}.\\
a) We start with a graph with only the source node:
    \begin{center}
    \includegraphics[scale=1]{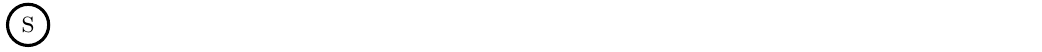}  	
    \end{center}
b) At the first iteration, we add arcs associated with the first
item, which has weight~5. There is space only for one item of this size.
    \begin{center}
    \includegraphics[scale=1]{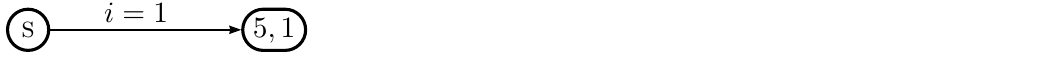}  	
    \end{center}
c) At the second iteration, we add arcs for the second item,
which has weight~3. The demand does not allow us to form paths
with more than~1 consecutive arc of this size.
    \begin{center}
    \includegraphics[scale=1]{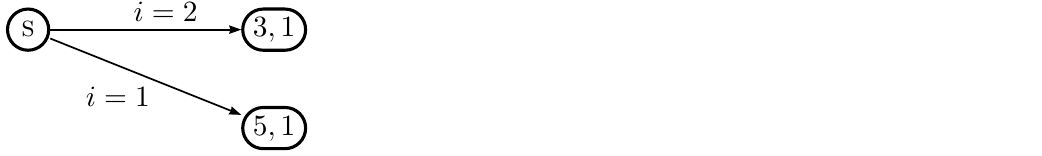}
    \end{center}
d) At the third iteration, we add the arcs associated
with the third item, which has weight~2. Since the demand of this item
is~2, we can add paths (starting from previously existing nodes) with at most~2 items of size~2.
    \begin{center}
    \includegraphics[scale=1]{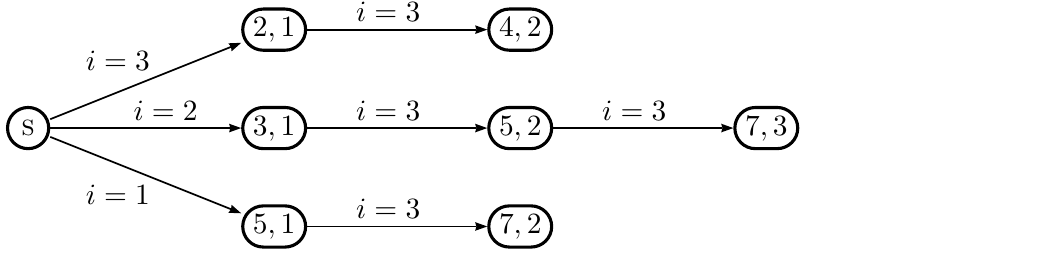}  
    \end{center}
e) Finally, we add the loss arcs connecting each node,
except the source, to the target.
The resulting graph corresponds to the model to be solved by a general-purpose
mixed-integer optimization solver.
    \begin{center}
    \includegraphics[scale=1]{cbpgraph.pdf}    
    \end{center}
\end{minipage}}
\end{figure}

The source vertex (\vS) is labeled with~$0$ in every dimension.
Arcs associated with items of weight $(w_i^1,w_i^2,\ldots,w_i^p)$ 
are created between vertices $(a^1,a^2,\ldots,a^p)$ and 
$(a^1+w_i^1,a^2+w_i^2,\ldots,a^p+w_i^p)$.
Since we just need to consider paths that respect a fixed order,
we can have an arc with tail in a node only if it is either the source node
or the head of an arc associated with a previous item
(according to the order defined in Definition~\ref{def:order}).
Our algorithm to construct the graph relies on this rule.
Initially, there is only the source node.
For each item, we insert in the graph arcs associated with the item
starting from all previously existing nodes .
After processing an item, we add to the graph the set
of new nodes that appeared as heads of new arcs.
This process is repeated for every item and
in the end we just need to connect every node, except the source,
to the target.
Using this algorihtm, the graph can be constructed in 
pseudo-polynomial time~$O(|V|m)$, where~$|V|$ is the number of
vertices in the graph.
Figure~\ref{fig:graphalg1} shows a small example of the construction of a graph
using this method.
More details on algorithms for graph construction 
are given in \cite{MThesisBrandao}.

\FloatBarrier
\subsection{Breaking symmetry}
\label{sec:breaksymmetry}

\begin{figure}[h!tbp]
\caption{Initial graph/Step-1 graph (with symmetry).\label{ex:compression0}}

  \centering
  \includegraphics[scale=1]{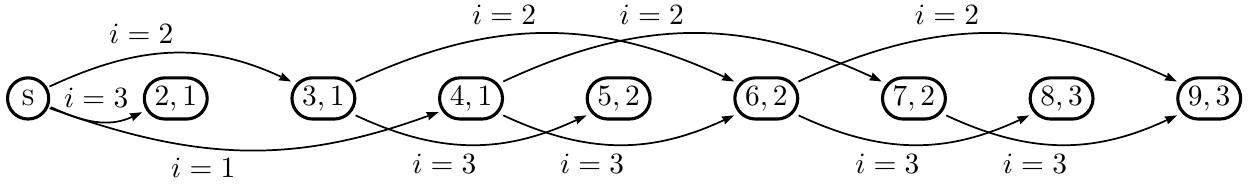}

\caption*{\footnotesize 
Graph corresponding to a cardinality constrained bin packing instance 
with bins of capacity~$W = 9$, cardinality limit~3 and items of sizes 
4, 3, 2 with demands 1, 3, 1, respectively.
}
\end{figure}

In Section~\ref{sec:compression}, we will present a three-step graph 
compression method whose first step consists of breaking the symmetry.
Let us consider a cardinality constrained bin packing instance with bins of capacity~$W = 9$,
cardinality limit~3 and items of sizes 4, 3, 2 with demands 1, 3, 1, respectively.
Figure~\ref{ex:compression0} shows the Step-1 graph produced by the 
graph construction algorithm without the final loss arcs.
This graph contains symmetry.
For instance, the paths (\vS,(4,1),i=1) ((4,1),(7,2),i=2) ((7,2),(9,3),i=3) and 
(s,(4,1),i=1) ((4,1),(6,2),i=3) ((6,2),(9,3),i=2) correspond to the same
pattern with one item of each size, but the second one does not respect 
the order defined in Definition~\ref{def:order}.

An easy way to break symmetry is to divide the graph into levels,
one level for each different item. 
We introduce in each node a new dimension that indicates the level
where it belongs.
For example, for the two-dimensional case, 
nodes $(a^{\prime}, b^{\prime})$ are transformed into sets of nodes 
$\{(a^{\prime}, b^{\prime},i^{\prime}),$ $(a^{\prime}, b^{\prime},i^{\prime\prime}),$ $\ldots\}$.
Each set has at most one node per level; 
nodes in consecutive levels are connected by loss arcs.
Arcs $((a^{\prime}, b^{\prime}),(a^{\prime\prime},b^{\prime\prime}),i)$ are
transformed into arcs $((a^{\prime}, b^{\prime},i),(a^{\prime\prime},b^{\prime\prime},i),i)$.
In level $i$, we have only arcs associated with items of weight $w_i$.
If we connect a node $(a^{\prime}, b^{\prime},i^{\prime})$ to a node $(a^{\prime}, b^{\prime},i^{\prime\prime})$
only in case $i^{\prime} < i^{\prime\prime}$, we ensure that every path will respect the 
order (defined in Definition~\ref{def:order}) and thus there is no symmetry.
Recall that the initial graph must contain every valid packing pattern
(respecting the order) represented as a path from~$\vS$ to~$\vT$.

\begin{figure}[h!tbp]
\caption{Graph with levels/Step-2 graph (without symmetry).\label{ex:compression1}}

  \centering
  \includegraphics[scale=1]{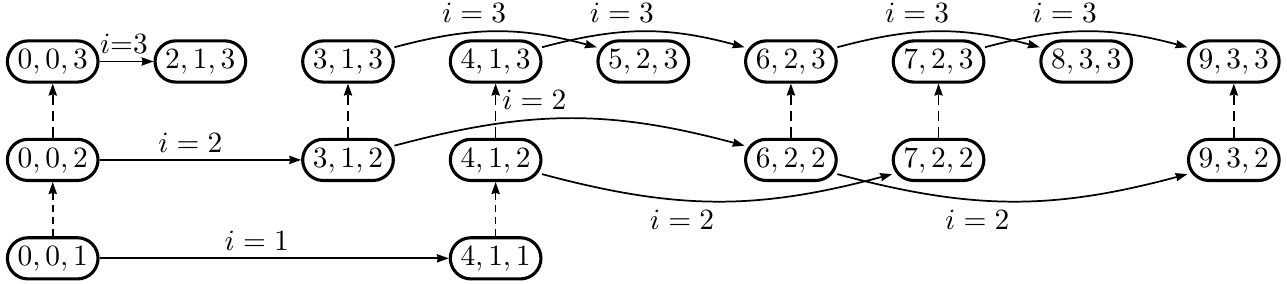}

\caption*{\footnotesize All the patterns respect the 
order since there are no arcs from higher levels
to lower levels. Moreover, it is also easy to check that no valid pattern
was removed.
In this graph, we consider $\vS = (0,0,1)$
since it is the only node without arcs incident to it.
}
\end{figure}

Figure~\ref{ex:compression1} shows the graph with levels (Step-2 graph)
that results from applying this symmetry breaking method to the
graph in Figure~\ref{ex:compression0}.
Although there is no symmetry, there are still
patterns that use some items more than their demand.
To avoid this, other alternatives to break symmetry could be used;
however this method is appropriate for the sake of simplicity
and speed.

\begin{property}
No valid pattern (respecting the order) is removed by
breaking symmetry with levels as long as the original graph
contains every valid packing pattern (respecting the order)
represented as a path from~$\vS$ to~$\vT$.
\end{property}

\begin{proof} For every valid packing pattern (respecting the order)
in the initial graph, there is a path in the Step-2 graph corresponding
to the same pattern. Note that every path can be seen as a sequence
of consecutive arcs. Let $(v_1, v_2, i_1)$ and $(v_2, v_3, i_2)$
be a pair of consecutive arcs in any valid packing pattern
in the Step-1 graph.
In the Step-2 graph, these arcs appear as $((v_1, i_1), (v_2, i_1), i_1)$
and $((v_2, i_2), (v_3, i_2), i_2)$. If $i_1 = i_2$, the pair of consecutive
arcs appear connected at the same level. 
If not, and given that $(v_2, i_1)$ and $(v_2, i_2)$ were created from~$v_2$,
a set of loss arcs
that connect nodes in different levels 
exists between them and hence there is again a sequence of arcs for that part of the pattern.
\end{proof}

\FloatBarrier
\subsection{Graph compression}
\label{sec:compression}

Symmetry may be helpful as long as it leads to large reductions
in the graph size. In this section we show how to reduce the graph size
by taking advantage of common sub-patterns that can be 
represented by a single sub-graph.
This method may increase symmetry, but it usually helps
by reducing dramatically the graph size.
The graph compression method is composed of three steps,
the first of which was presented in Section~\ref{sec:breaksymmetry}.

In the graphs we have seen so far, a node label $(a^1, a^2, \ldots, a^p)$
means that, for every dimension~$d$, every sub-pattern from the source to the node
uses~$a^d$ space in that dimension.
This means that~$a^d$ corresponds to the length of the longest
path from the source to the node in dimension $d$.
Similarly, the longest path to the target can also be used
as a label and nodes with the same label can be combined into one single node.
In the main compression step, 
a new graph is constructed using the longest path to the target in each
dimension.
This usually allows large reductions in the graph size.
This reduction can be improved by breaking symmetry first (as described
in the previous section),
which allows us to consider only paths to the target with a specific order.

\begin{figure}[h!tbp]
\caption{Step-3 graph (after the main compression step).\label{ex:compression2}}

  \centering
  \includegraphics[scale=1]{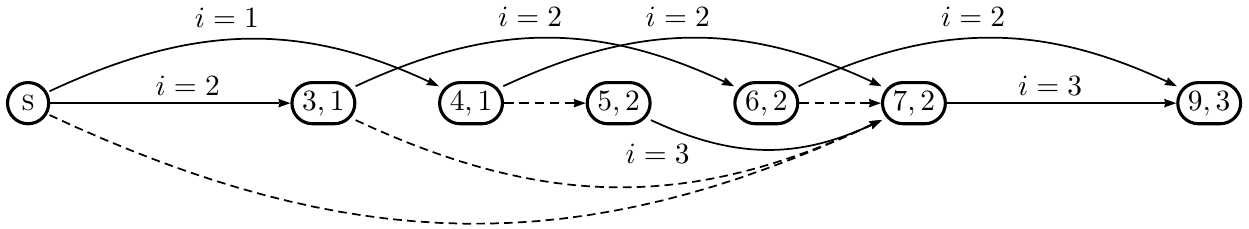}

\caption*{\footnotesize The Step-3 graph has 8 nodes and 17 arcs 
(considering also the final loss arcs connecting internal nodes to~$\vT$).
}
\end{figure}

The main compression step is applied to the Step-2 graph.
In the Step-3 graph, the longest paths to the target in each dimension is used
to relabel the nodes, dropping the level dimension of each node.
Let ($\varphi^1(u)$, $\varphi^2(u)$, \ldots, $\varphi^p(u)$) be the new label of 
node~$u$ in the Step-3 graph, where
\begin{alignat}{3}
\varphi^d(u) & = & \left\{ \begin{array}{ll}
                0 & \mbox{if }  u = \vS, \\
                W^d & \mbox{if }  u = \vT, \\
                \min_{(u',v,i) \in A:u'=u}\{\varphi^d(v) - w_i^d\}  & \mbox{otherwise.}\\
                \end{array}\right.           
\end{alignat}
For the sake of simplicity, we define $w_0^d$ for loss arcs as zero in every dimension, 
$w_0^0=w_0^1=\ldots=w_0^p=0$.
In the paths from $\vS$ to $\vT$ in the Step-2 graph usually there is float in some dimension.
In this process, we are moving this float as much as possible to the beginning of the path.
The label in each dimension of every node~$u$ (except~$\vS$) corresponds to the highest position
where the sub-patterns from~$u$ to~$\vT$ can start in each dimension
so that capacity constraints are not violated.
By using these labels we are allowing arcs to be longer than the items
to which they are associated.
We use dynamic programming to compute these labels in linear time.
Figure~\ref{ex:compression2} shows the Step-3 graph that
results from applying the main compression step to the graph of Figure~\ref{ex:compression1}.
Even in this small instance, a few nodes and arcs were removed
comparing with the initial graph of Figure~\ref{ex:compression0}.

\begin{figure}[h!tbp]
\caption{Step-4 graph (after the last compression step).\label{ex:compression3}}

  \centering
  \includegraphics[scale=1]{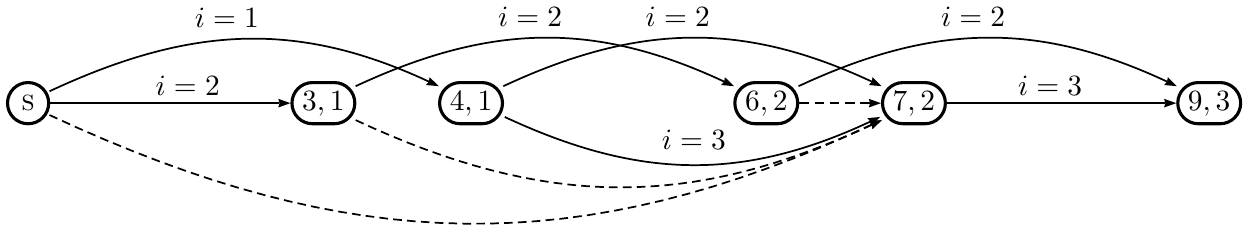}

\caption*{\footnotesize The Step-4 graph has 7 nodes and 15 arcs
(considering also the final loss arcs connecting internal nodes to~$\vT$).
In this case, the only difference from the Step-3 graph is the
node~$(5,2)$ that collapsed with the node~$(4,1)$.
The initial Step-1 graph had 9 nodes and 18 arcs.}
\end{figure}

Finally, in the last compression step, a new graph is constructed once more. 
In order to try to reduce the graph size even more, 
we relabel the graph once more using the longest paths from the source in each dimension.
Let ($\psi^1(v)$, $\psi^2(v)$, \ldots, $\psi^p(v)$) be the label of node~$v$ in the
Step-4 graph, where
\begin{alignat}{3}
\psi^d(v) & = & \left\{ \begin{array}{ll}
                0 & \mbox{if }  v = \vS, \\
                \max_{(u,v',i) \in A:v'=v}\{\psi^{d}(u) + w_i^d\}  & \mbox{otherwise.}\\
                \end{array}\right.           
\end{alignat}
Figure~\ref{ex:compression3} shows the Step-4 graph
without the final loss arcs. This last compression step is not 
as important as the main compression step, but it is
easy to compute and usually removes many nodes and arcs.

Note that, in this case, the initial Step-1 graph had some symmetry
and the final Step-4 graph does not contain any symmetry.
Graph compression may increase symmetry in some situations,
but in practice this is not a problem as long as it leads to large 
reductions in the graph size.
Since we are dealing with a very small instance,
the improvement obtained by compression is not as substantial as in large instances. 
For instance, in the standard BPP instance HARD4 from \cite{Scholl1997627} 
with 200 items of 198 different sizes
and bins of capacity $100,000$, we obtained reductions of $97\%$ in the number of vertices
and $95\%$ in the number of the arcs. The resulting model was solved in a few seconds.
Without graph compression it would be much harder
to solve this kind of instances using the arc-flow formulation.

\begin{property}[Graph compression 1]
Non-redundant patterns of the initial graph are not removed 
by graph compression.
\end{property}
\begin{proof}
Any pattern in the initial graph is represented by a path
from~$\vS$ to~$\vT$. Consider a path $\pi = v_1v_2\ldots v_n$.
Graph compression will reduce the graph size by relabeling nodes,
collapsing nodes that receive the same label
and removing self-loops (loss arcs).
Therefore, the path $\pi^{\prime} = \phi(v_1)\phi(v_2)\ldots\phi(v_n)$,
where~$\phi$ is the map between the initial and final labels, 
will represent exactly the same pattern as the one represented by the path~$\pi$,
possibly with some self-loops that do not affect the patterns.
\end{proof}

\begin{property}[Graph compression 2]
Graph compression will not introduce any invalid pattern.
\end{property}
\begin{proof}
An invalid pattern consists of a set of items whose sum of weights 
exceeds the bin capacity in some dimension.
A pattern is formed from a path in the graph
and its total weight in each dimension is the sum of weights of the items in the path.
The main compression step just relabels every node~$u$ (except~$\vS$) in each
dimension with the highest position
where the sub-patterns from~$u$ to~$\vT$ can start
so that capacity constraints are not violated.
The node~$\vT$ is labeled with $(W^1, W^2, \ldots, W^p)$
and no label will be smaller than $(0, 0, \ldots, 0)$
since all the patterns in the input graph are required to be valid.
However, we could have invalid patterns, even with all the nodes
labeled between zero and $(W^1, W^2, \ldots, W^p)$, if
an arc had length smaller than the item it represents,
but this is not possible. The label of every internal node $u$
is given by $(\varphi^1(u), \varphi^2(u), \ldots, \varphi^p(u))$,
where $\varphi^d(u) = \min\{\varphi^d(v) - w_i^d\ |\ (u',v,i) \in A,\ u'=u\}$;
for every node~$v$ such that there is an arc between~$u$ and~$v$,
the difference between their labels
is at least the weight in each dimension of the item associated with the arc.
Therefore, no invalid patterns are introduced.
An analogous proof can be derived for the last compression step.
\end{proof}

\FloatBarrier
\subsection{Building Step-3 graphs directly}

As we said in Section~\ref{sec:multidimgraphs}, 
in the $p$-dimensional case, the number of vertices
and arcs is limited by $O(\prod_{d=1}^{p} (W^d+1))$ 
and $O(m\prod_{d=1}^{p} (W^d+1))$, respectively.
Our graph compression method usually leads to very high compression
ratios and in some cases it may lead to final graphs
hundreds of times smaller than the initial ones.
The size of the initial graph can be the limiting factor
and hence its construction should be avoided.

In practice, we build the Step-3 graph directly
in order to avoid the construction of huge Step-1 and Step-2 graphs
that may have many millions of vertices and arcs.
Algorithm~\ref{alg:build} uses dynamic programming to build the Step-3 graph
recursively over the structure of the Step-2 graph
(without building it).
The base idea for this algorithm comes from the fact that
in the main compression step the label of any node only depends
on the labels of the two nodes to which it is connected 
(a node in its level and another in the level above).
After directly building the Step-3 graph from the instance's 
data using this algorithm, we just need 
to apply the last compression step to obtain the final graph. 
In practice, this method allows obtaining 
arc-flow models even for large 
benchmark instances quickly.

The dynamic programming states are identified by
the space used in each dimension~($x^d$, for $d=1,\ldots,p$), the current item~($i$)
and the number of times it was already used~($c$).
In order to reduce the number of states, we lift
each state by solving (using dynamic programming) knapsack/longest-path problems 
in each dimension considering the remaining items;
we try to increase the space used in each dimension
to its highest value considering the valid packing patterns for the remaining items.
By solving multi-constraint knapsack problems in each dimension,
we would obtain directly the corresponding label in the Step-3 graph;
however, it is very expensive to solve this
problem many times.
By lifting states as a result of solving one-dimensional knapsack problems,
we obtain a good approximation in a reasonable amount of time, 
and it usually leads to a substantial reduction in the number of states.
\begin{algorithm}[!h]
\caption{Direct Step-3 Graph Construction Algorithm}\label{alg:build}
\SetKwInOut{Input}{input}
\SetKwInOut{Output}{output}

\SetKwFunction{Build}{build}
\SetKwFunction{BuildGraph}{buildGraph}
\SetKwFunction{Sort}{sort}
\SetKwFunction{Key}{key}
\SetKwFunction{Reversed}{reversed}
\SetKwFunction{Lift}{lift}
\SetKwFunction{HPos}{highestPosition}
\SetKwData{vardp}{dp}
\SetKwData{varx}{x}
\SetKwData{vari}{i}
\SetKwData{varmx}{mx}
\SetKwData{varmi}{mi}
\SetKwData{NIL}{NIL}

\SetInd{0.5em}{1.5em}

\SetKwBlock{Function}{}{}

\Input{$m$ - number of different items; $w$ - item sizes; $b$ - demand; $W$ - capacity}
\Output{$V$ - set of vertices; $A$ - set of arcs; $\vS$ - source; $\vT$ - target}

\textbf{function} $\BuildGraph(m, w, b, W)$:
\Function{
$\vardp[x', i', c'] \gets \NIL, \mbox{ \bf for all } x', i', c'$\tcp*[r]{dynamic programming table}
\textbf{function} $\Lift(x, i, c)$: \tcp*[f]{auxiliary function: lift \vardp states solving knapsack/longest-path problems in each dimension}
\Function{
    \textbf{function} $\HPos(d, x, i, c)$:
    \Function{
       \Return \\ \vspace{-4.5mm} \hspace{12mm}
       $
        \begin{array}{lll}
        \mbox{min} & W^d-\sum_{j = i}^{m} w_j^d y_j \\  
        \mbox{s.t.} & \sum_{j = i}^{m} w_j^d y_j \leq W^d-x^d,\\
                     & y_i \leq b_i-c,  &  \\
                     & y_j \leq b_j, \qquad j=i+1,...,m &  \\
                     & y_j \geq 0, \mbox{ integer},\  j=i,...,m;   &
        \end{array}
        $
    }
    \Return $(\HPos(1, x, i, c),\ldots, \HPos(p, x, i, c))$\;    
}
$V \gets \{\ \}$; $A \gets \{\ \}$\;
\textbf{function} $\Build(x, i, c)$:
\Function{
$x \gets \Lift(x,i,c)$\tcp*[r]{lift $x$ in order to reduce the number of \vardp states}\label{alg:lift}
\If(\tcp*[f]{avoid repeating work }){$\vardp[x,i,c] \neq \NIL$}{
    \Return $\vardp[x,i,c]$\;
}
$u \gets (W^1,...,W^p)$\;
\If(\tcp*[f]{option 1: do not use the current item (go to the level above)}){$i < m$}{
    $up_x \gets \Build(x, i+1, 0)$\;
    $u \gets up_x$\;
}
\If(\tcp*[f]{option 2: use the current item}){$c < b_i \mbox{ \bf and }  x^d+w_i^d \leq W^d \mbox{ \bf for all } 1 \leq d \leq m$}{    
    $v \gets \Build((x^1+w_i^1,\ldots,x^p+w_i^p), i+1, c+1)$\;
    $u \gets (\min(u^1, v^1-w_i^1),\ldots,\min(u^p, v^p-w_i^p))$\;
    $A \gets A \cup \{(u, v,i)\}$\;
    $V \gets V \cup \{u, v\}$\;
    \If(){$i < m${ \bf and }$u \neq up$}{
        $A \gets A \cup \{(u,up_x,0)\}$\;
        $V \gets V \cup \{up_x\}$\;        
    }
}
$\vardp[x,i,c] \gets u$\;

\Return $u$\;
}
$\vS \gets \Build(x=(0,\ldots,0),i=0,c=0)$\tcp*{build the graph}
$V \gets V \cup \{\vT\}$\;
$A \gets A \cup \{(u, \vT, 0)\ |\ u \in V\setminus \{\vS, \vT\} \}$\tcp*{connect internal nodes to the target}
\Return $(G = (V, A), \vS, \vT)$\;
}
\end{algorithm}

\subsection{Integrality gap and heuristic solutions}
\label{sec:gap}

There have been many studies 
(see, e.g., \citealt{Scheithauer95themodified}, \citealt{Scheithauer1997}) 
about the integrality gap for the BPP, many of them about
the integrality gap using Gilmore-Gomory's model.
Our proposed arc-flow formulation
is equivalent to  Gilmore-Gomory's model and hence
the lower bounds are the same when the same set of patterns is considered.
Therefore, the results found in these studies, which we summarize next, are also valid 
for the arc-flow formulation.

\begin{definition}[Integer Property]\label{def:IP}
A linear integer optimization problem P has the integer property (IP) if
\[
z_{ip}^*(E) = z_{lp}^*(E) \mbox{\ for every instance\ } E \in P
\]
\end{definition}

\begin{definition}[Integer Round-Up Property]\label{def:IRUP}
A linear integer optimization problem P has the integer round-up property (IRUP) if
\[
z_{ip}^*(E) = \lceil z_{lp}^*(E) \rceil \mbox{\ for every instance\ } E \in P
\]
\end{definition}

\begin{definition}[Modified Integer Round-Up Property]\label{def:MIRUP}
A linear integer optimization problem P has the modified integer round-up property (MIRUP) if
\[
z_{ip}^*(E) = \lceil z_{lp}^*(E) \rceil + 1 \mbox{\ for every instance\ } E \in P
\]
\end{definition}

\cite{Rietz2002_non-IRUP} describe families of instances 
of the one-dimensional cutting stock problem
without the integer round-up property. 
One of the families is the
so-called divisible case,
which was firstly proposed by \cite{nica_counterexample}, 
where every item size~$w_i$ is a factor of the bin capacity~$W$.
Since the method we propose
usually solves bin packing problems quickly, it was used to solve
millions of instances from this family
keeping track of the largest gap found, which was~$1.0378\ldots$, 
the same gap as the one found by
\cite{Scheithauer1997}.

\cite{gau_gap} presents an instance with a gap of~$1.0666$.
The largest gap known so far is~$7/6$ and it was found by~\cite{Rietz2002_bounds}.
\cite{Scheithauer1997} conjecture that
the general one-dimensional cutting stock problem has the modified
integer round-up property (MIRUP).  Moreover, instances for BPP and CSP
usually have the integer round-up
property (IRUP).  Concerning the results obtained using the arc-flow
formulation in several different cutting and packing problems
(see Section~\ref{sec:results}), most of the instances have the
IRUP, and no instance violated the MIRUP.  
The largest gap we found in all the instances 
from all the benchmark test data sets 
(except the case of the previous paragraph) was~1.0027.

The lower bound provided by the linear relaxation of the proposed arc-flow formulation is 
usually very tight for every problem we considered; hence, the
branch-and-bound process usually finds the optimal solution quickly (see Section~\ref{sec:results}).
In fact, in our experiments, the large majority of models were solved at
the root node.
Moreover, very good solutions are usually found when rounding
the linear programming (LP) solution.
Rounding up the fractional variables of the LP
solution of Gilmore and Gomory's model
guarantees a heuristic solution of value at most~$z_{lp}^*+m$, where~$z_{lp}^*$
is the optimum value of the linear relaxation,
since we need to round up at most one variable for each different item size
in order to obtain a valid integer solution.
Rounding up fractional flow paths
of the linear relaxation of the arc-flow formulation
gives the same guarantee.
\cite{springerlink:10.1007/BF01539705} present
more elaborate rounding heuristics for Gilmore and Gomory's model
that usually lead to the optimal solution
in cutting stock instances. Note that these rounding
heuristics usually work well in cutting stock instances where the demands
are large, but they may have a poor performance in bin packing instances
where the values of variables are often a fraction of unity.
In our model, we try to overcome the problems introduced by low demands
with upper bounds on the variable values (constraints (\ref{eq:new5}) of the
general arc-flow model), and by requiring the demand of items with demand 
one to be satisfied exactly (constraints (\ref{eq:new4})).

\FloatBarrier
\section{Applications and results}
\label{sec:results}

In this section, we present the results obtained using the arc-flow formulation
in several cutting and packing problems.
CPU times were obtained using a computer with two Quad-Core Intel Xeon at
2.66GHz, running Mac OS X 10.8.0, with 16 GBytes of memory
(though only  a few of the hardest instances required more than 4 GBytes of memory). 
The graph construction algorithm was implemented in~\texttt{C++}, 
and the resulting MIP was solved using \texttt{Gurobi}~5.0.0, a state-of-the-art
mixed integer programming solver. 
The parameters used in \texttt{Gurobi} were Threads~=~1 (single
thread), Presolve~=~1 (conservative), Method~=~2 (Interior point methods),
MIPFocus~=~1 (feasible solutions), Heuristics~=~1, MIPGap~=~0, MIPGapAbs~=~$1-10^{-5}$ 
and the remaining parameters were \texttt{Gurobi}'s default values. The
branch-and-cut solver used in \texttt{Gurobi} uses a series of cuts; in our
models, the most frequently used were Gomory, Zero half and MIR.
The source code is available online\footnote{\url{http://www.dcc.fc.up.pt/\~fdabrandao/code}}.
For more detailed results please consult the web-page~\cite{FlowResults}, 
which contains all the instances, results and log files.

\FloatBarrier
\subsection{$p$-dimensional vector packing}

In the $p$-dimensional vector packing problem, bins and items have $p$~independent
dimensions.
For each dimension~$d$, let~$w_i^d$ be the weight of item~$i$ and~$W^d$ the bin capacity.
The set $S$ of valid patterns for this problem is defined as follows:
\begin{equation}
  A = \begin{bmatrix}
    w_{1}^{1} & \ldots & w_{m}^{1}  \\
    \svdots & & \svdots\\
    w_{1}^{p} & \ldots & w_{m}^{p} \\
  \end{bmatrix}
  \qquad
  L = \begin{bmatrix}
    W^{1}\\
    \svdots \\
    W^{p}\\
  \end{bmatrix}
  \qquad 
  S = \{\XVEC \in \mathbb{N}_{0}^{m} : A\XVEC \leq L\}
\end{equation}
where~$A$ is the matrix of weights and $L$ is the vector of capacities.
The set $S$ is the set of valid packing patterns that satisfy 
all the following knapsack constraints:
\begin{alignat}{6}
w_1^1 x_1 &+& w_2^1 x_2 &+& \ldots &+& w_m^1 x_m &\leq& W^1 &&\label{vpconbegin} \\
w_1^2 x_1 &+& w_2^2 x_2 &+& \ldots &+& w_m^2 x_m &\leq& W^2 && \\
\centermathcell{\svdots} && \centermathcell{\svdots} &&  && \centermathcell{\svdots} && \centermathcell{\svdots} \\
w_1^p x_1 &+& w_2^p x_2 &+& \ldots &+& w_m^p x_m &\leq& W^p &&\label{vpconend1}\\
\omit\rlap{$x_i \geq 0,\ \mbox{integer},\ i=1,\ldots,m,$}\label{vpconend2} 
\end{alignat}
By replacing the subproblem (\ref{eq:cuttingstock4})
of model (\ref{eq:cuttingstock1})-(\ref{eq:cuttingstock3})
by the knapsack constraints (\ref{vpconbegin})-(\ref{vpconend2}), we obtain 
a variant of Gilmore-Gomory's model for $p$-dimensional vector packing.
In our method, these patterns are derived from paths in a graph.
We start by building a graph~$G=(V,A)$ containing every valid packing pattern represented 
as a path from the source to the target, 
and an optimal linear combination of patterns is obtained
through the solution of the general arc-flow formulation (\ref{eq:new1})-(\ref{eq:new6}) over $G$.

Two-constraint bin packing (2CBP) is a bi-dimensional vector packing problem.
We used the proposed arc-flow formulation to solve 330 of the 400 instances from
the \cite{TCBPInstances} two-constraint bin packing test data set, which was
proposed by \cite{Caprara:2001:LBA:508252.508254}.  Table~\ref{tab:2dvbp}
summarizes the results.  This data set has several sizes for each class, 
each pair (class, size) having~10 instances.

\begin{table}[h!tbp]
  \centering\begin{threeparttable}[b]
    \scriptsize 
    \caption{Results for 2-dimensional vector packing.\label{tab:2dvbp}}
    \begin{tabular}{rrrrrrrrrrrrrrrrr}
      \hline
      \up\down
      & & \multicolumn{3}{c}{$n \in \{24,25\}$} & & \multicolumn{3}{c}{$n \in \{50,51\}$} & &
      \multicolumn{3}{c}{$n \in \{99,100\}$} & & \multicolumn{3}{c}{$n \in \{200,201\}$} \\
      \cline{3-5} \cline{7-9} \cline{11-13} \cline{15-17}
      \up\down
      class & & $\nbb$ & $\ttot$ & \#op &   & $\nbb$ & $\ttot$ & \#op &  & $\nbb$ & $\ttot$
      &  \#op &  & $\nbb$ & $\ttot$ &  \#op \\
      \hline
      \up
        1 &  & 0.0 & 0.12 & 0 &  & 0.0 & 1.62 & 0 &  & 0.0 & 66.96 & 5 &  & 0.0 & 7,601.35 & 7\\
        2 &  & 0.0 & 0.01 & 10 &  & 0.0 & 0.04 & 10 &  & 0.0 & 0.21 & 10 &  & 0.0 & 6.93 & 10\\
        3 &  & 0.0 & 0.01 & 10 &  & 0.0 & 0.02 & 10 &  & 0.0 & 0.05 & 10 &  & 0.0 & 0.20 & 10\\
        4 &  & 0.0 & 21.97 & 10 &  & - & - & - &  & - & - & - &  & - & - & -\\
        5 &  & 0.0 & 10.62 & 10 &  & - & - & - &  & - & - & - &  & - & - & -\\
        6 &  & 0.0 & 0.02 & 0 &  & 0.0 & 0.06 & 1 &  & 0.0 & 0.31 & 5 &  & 0.0 & 4.75 & 8\\
        7 &  & 0.0 & 0.03 & 0 &  & 0.0 & 0.14 & 1 &  & 0.3 & 1.69 & 7 &  & 31.4 & 14.01 & 3\\
        8 &  & 0.0 & 0.01 & 10 &  & 0.0 & 0.02 & 10 &  & 0.0 & 0.07 & 10 &  & 0.0 & 0.24 & 10\\
        9 &  & 0.0 & 0.10 & 0 &  & 0.0 & 0.66 & 1 &  & 0.0 & 28.11 & 10 &  & - & - & -\\
        \down
        10 &  & 0.0 & 0.02 & 0 &  & 0.0 & 0.10 & 0 &  & 0.0 & 0.66 & 0 &  & 226.9 & 155.43 & 0\\
      \hline
    \end{tabular}
    \caption*{\footnotesize 
        \footnotesize 
        $n$~-~total number of items;
        $\nbb$~-~average number of nodes explored in the branch-and-bound procedure; 
        $\ttot$~-~average run time in seconds;
        \#op~-~number of previously open instances solved.
    }    
  \end{threeparttable}
\end{table}

Note that among the instances presented in Table~\ref{tab:2dvbp} there were 188
instances with no previously known optimum.  The arc-flow formulation
allowed the solution of 330 instances out of the 400 instances; hence, 
there remain 70 open instances.
The graph compression algorithm is remarkably effective on all of this
subset of two-constraint bin packing instances. In many of the instances, graph compression
allowed the remotion of more than $90\%$ of the vertices and arcs;
without it, it would not be viable to
solve many of these instances within a reasonable amount of time. 

One may ask why so many of these instances could not be solved before.  A
reasonable explanation may be the fact that, for example, in instances from
classes 2, 3 and 8, the lower bound provided by the linear relaxation of
assignment-based formulations is rather loose. 
Also, in some of these instances, the average number of items per
bin in the optimal solution is reasonably large, which leads to an extremely
large number of possible patterns. 
Results presented in \cite{Caprara:2001:LBA:508252.508254} show
that solving the linear relaxation of many of these instances using
Gilmore-Gomory's model and column-generation was not possible within 100,000 seconds (more than
27 hours).  The computer they used is much slower than the one we
used here, but this shows how hard it can be to compute a linear
relaxation of Gilmore-Gomory's model as the number of patterns increases and the
subproblems become harder to solve.  With the arc-flow formulation, graph compression
leads to very large reductions in the graph size and allows us to represent all
these patterns in reasonably small graphs.  None of the instances from these
classes has been solved before.

The classes 4 and 5 are hard even for our arc-flow formulation, 
because of the large number of items that fit in a single bin. 
For this type of instances, assignment-based formulations tend to perform better 
due to the heuristics used by integer programming solvers, since
it is usually easier to find very good solutions when the number of items
that fit in each bin is large.
Most of the instances that remain open belong to these two classes.  
The 7 subclasses that we did not solve contain at
least one instance that takes more than 12 hours to be solved exactly.
The average run time in the 330 solved instances was 4 minutes
and none of these instances took longer than 5 hours to be solved exactly.

In order to test the behavior of the arc-flow formulation
in instances with the same characteristics and more dimensions,
we created 20-dimensional vector packing instances by combining the ten
2-dimensional vector packing instances of each subclass (class, size)
into one instance.
Table~\ref{tab:20dvbp} summarizes the results.
The arc-flow formulation allowed the solution of 33 out of the 40 
instances. 
The same subclasses that were solved in the 2-dimensional case were
also solved in the 20-dimensional case.
Moreover, some 20-dimensional instances were easier to solve
than the original 2-dimensional ones due to the reduction
in the number of valid packing patterns.
The 7 instances that were not solved within a 12 hour time limit
are mostly instances in which the patterns are very long.
The average run time in the other 33 solved instances was 48 seconds,
and none of these instances took longer than 23 minutes to be solved exactly.

\begin{table}[h!tbp]
  \centering\begin{threeparttable}[b]
    \scriptsize 
    \caption{Results for 20-dimensional vector packing.\label{tab:20dvbp}}
    \begin{tabular}{rrrrrrrrrrrrr}
      \hline
      \up\down
      & & \multicolumn{2}{c}{$n \in \{24,25\}$} & & \multicolumn{2}{c}{$n \in \{50,51\}$} & &
      \multicolumn{2}{c}{$n \in \{99,100\}$} & & \multicolumn{2}{c}{$n \in \{200,201\}$} \\
      \cline{3-4} \cline{6-7} \cline{9-10} \cline{12-13}
      \up\down
      class & & $\nbb$ & $\ttot$ &   & $\nbb$ & $\ttot$ &  & $\nbb$ & $\ttot$
      &  & $\nbb$ & $\ttot$ \\
      \hline
      \up
        1 &  & 0 & 0.09 &  & 0 & 0.81 &  & 0 & 36.28 &  & 0 & 1,374.18\\
        2 &  & 0 & 0.01 &  & 0 & 0.01 &  & 0 & 0.01 &  & 0 & 0.01\\
        3 &  & 0 & 0.01 &  & 0 & 0.01 &  & 0 & 0.01 &  & 0 & 0.01\\
        4 &  & 0 & 50.27 &  & - & - &  & - & - &  & - & -\\
        5 &  & 0 & 73.20 &  & - & - &  & - & - &  & - & -\\
        6 &  & 0 & 0.01 &  & 0 & 0.02 &  & 0 & 0.05 &  & 0 & 0.19\\
        7 &  & 0 & 0.01 &  & 0 & 0.03 &  & 0 & 0.06 &  & 0 & 0.19\\
        8 &  & 0 & 0.01 &  & 0 & 0.01 &  & 0 & 0.03 &  & 0 & 0.10\\
        9 &  & 0 & 0.05 &  & 0 & 0.37 &  & 0 & 12.80 &  & - & -\\
        \down
        10 &  & 0 & 0.02 &  & 0 & 0.11 &  & 0 & 0.91 &  & 0 & 14.52\\
      \hline
    \end{tabular}
    \caption*{\footnotesize 
        $n$~-~total number of items;        
        $\nbb$~-~number of nodes explored in the branch-and-bound procedure; 
        $\ttot$~-~run time in seconds.
    }
  \end{threeparttable}
\end{table}

\FloatBarrier
\subsection{Graph coloring}
\label{sec:color}

The graph coloring problem is a combinatorial NP-hard problem (see, e.g.,
~\citealt{Garey:1979:CIG:578533}) in which one has to assign a color to each vertex of a
graph in such way that no two adjacent vertices share the same color
using the minimum number of colors.

Graph coloring can be reduced to vector packing in several ways. 
Let variables $x_i$ of constraints (\ref{vpconbegin})-(\ref{vpconend1}) 
represent whether or not vertex $i$ appears in a given pattern 
(each pattern corresponds to a set of vertices that can share the same color).
Considering each color as a bin and each vertex as an item with demand one, the following
reductions are valid:
\begin{itemize}
\item Adjacenty constraints: For each pair of adjacent vertices $i$ and $j$, there is an adjacency
  constraint $x_i+x_j \leq 1$. Each adjacency constraint can be represented by a
  dimension $k$ of capacity $W^k=1$, with $w_i^k = w_j^k = 1$.
\item Degree constraints: Let $\Deg(i)$ and $\Adj(i)$ be the degree and the list of adjacent
  vertices of vertex $i$, respectively.
  For each vertex $i$, there is a constraint $\Deg(i) x_i + \sum_{j \in \Adj(i)}
  x_j \leq \Deg(i)$.  
  Each constraint can be represented by a dimension $k$ of
  capacity $W^k=\Deg(i)$, with $w_i^k = \Deg(i)$ and $w_j^k=1$ for every $j \in
  \Adj(i)$.
\item Clique constraints: For each clique $C$, there is a constraint $\sum_{i \in C} x_i \leq 1$.
  Each clique constraint can be represented by a dimension $k$ of capacity $W^k
  = 1$, 
  with $w_i^k=1$ for every $i \in C$. \cite{BronKerbosch1973}'s algorithm
  can be used to decompose the graph into maximal cliques.
\end{itemize}

Using any of the three reductions above, a vector packing solution with $z$ bins
corresponds to a graph coloring solution with $z$ colors.  Different reductions
result in different vector packing instances that can be harder or easier to
solve.  According to our experiments, it is usually a good idea to choose
reductions that lead to vector packing instances with fewer dimensions.
Figure~\ref{fig:graphcolorexample} illustrates the three
reductions defined above.

\begin{figure}[h!tbp]
\caption{Graph coloring reductions to vector packing.\label{fig:graphcolorexample}}
\centering
\fbox{
	\hspace{-5mm}
	\begin{minipage}[t]{0.28\linewidth}\centering
	Graph coloring instance:\\
    \vspace{2mm}
	\begin{tikzpicture}
	\node (a) at (0bp,20bp) [draw,circle,] {$x_1$};
	\node (b) at (0bp,-20bp) [draw,circle,] {$x_2$};
	\node (c) at (40bp,0bp) [draw,circle,] {$x_3$};
	\node (d) at (80bp,0bp) [draw,circle,] {$x_4$};
	\draw [] (a) -- (b);
	\draw [] (b) -- (c);
	\draw [] (c) -- (a);
	\draw [] (c) -- (d);
	\end{tikzpicture}
	\end{minipage}  
	\begin{minipage}[t]{0.24\linewidth}\centering
	Adjacency constraints:
	\vspace{-3mm}
	\begin{alignat*}{5}
	x_1 &+& x_2 && && &\leq& 1,\\
	x_1 && &+& x_3 && &\leq& 1,\\
	&& x_2 &+& x_3 && &\leq& 1,\\
	&& && x_3 &+& x_4 &\leq& 1,\\
	\omit\rlap{$x_i \in \{0,1\},\ i=1..4$}
	\end{alignat*}
	\end{minipage}  
	\begin{minipage}[t]{0.23\linewidth}\centering
	Degree constraints:
	\vspace{-3mm} 
	\begin{alignat*}{5}
	2x_1 &+& x_2 &+& x_3 && &\leq& 2,\\
	2x_2 &+& x_1 &+& x_3 && &\leq& 2,\\
	3x_3 &+& x_1 &+& x_2 &+& x_4 &\leq& 3,\\
	1x_4 &+& x_3 && && &\leq& 1,\\
	\omit\rlap{$x_i \in \{0,1\},\ i=1..4$}
	\end{alignat*}
	\end{minipage}  		
	\begin{minipage}[t]{0.23\linewidth}\centering	
	Clique constraints: 
	\vspace{-3mm}
	\begin{alignat*}{5}
	x_1 &+& x_2 &+& x_3 &&  &\leq& 1,\\
	 &&  && x_3 &+& x_4 &\leq& 1,\\
	\omit\rlap{$x_i \in \{0,1\},\ i=1..4$}
	\end{alignat*}
	\end{minipage} 
	\hspace{-5mm} 
}	
\end{figure}

Note that, in graph coloring,  the length of the patterns
are usually very long, for instance, when the graphs are sparse. 
However, there are problems that can be reduced to 
graph coloring problems with reasonably short patterns,
and thus it may be possible to solve them using the proposed arc-flow model.
One of these problems is timetabling (see Section~\ref{sec:timetabling}).

Table~\ref{tab:coloring} shows the results for a small subset of graph
coloring instances from~\cite{ORLibrary}.
These instances correspond to queen graphs;
given a $q\times q$~chessboard, there are $q^2$~nodes (one for each square of the board)
which are connected by an edge if the corresponding squares are in the same row, 
column, or diagonal.
If there is a solution with at most $q$ colors, 
then it is possible to place $q$ sets of $q$ queens on the board 
so that no two queens of the same set are in the same row, column, or diagonal.
In these instances, $q$ limits the length of the color patterns and hence
the arc-flow formulation can be applied with success for small values of $q$.
We reduced these instances to vector packing through degree constraints;
while the arc-flow formulation solved these instances quickly,
the assignmnet-based formulation allowed us to solve
only the smallest one under an one-minute time limit.

\begin{table}[h!tbp]
\centering\begin{threeparttable}[b]
  \scriptsize 
  \caption{Results for graph coloring.\label{tab:coloring}}
  \begin{tabular}{lrrrrrrrrrrrrr}
    \hline
    \up\down
    name & $n$ & $e$ & $d$ & $\zip$ & $\zlp$ & $\#v$ & $\#a$ & $\tpp$ & $\tlp$ & $\tip$ & $\nbb$ & $\ttot$\\
    \hline
    \up
	 queen5\_5 & 25 & 320 & 25 & 5 & 5.00 & 95 & 378 & 0.01 & 0.00 & 0.01 & 0 & 0.02\\
	 queen6\_6 & 36 & 580 & 36 & 7 & 7.00 & 367 & 1,502 & 0.04 & 0.02 & 0.03 & 0 & 0.09\\
	 queen7\_7 & 49 & 952 & 49 & 7 & 7.00 & 1,559 & 6,571 & 0.23 & 0.11 & 0.12 & 0 & 0.46\\
	 \down
	 queen8\_8 & 64 & 1,456 & 64 & 9 & 8.44 & 7,799 & 34,280 & 1.69 & 1.78 & 5.41 & 0 & 8.88\\
    \hline
  \end{tabular}
  \caption*{\footnotesize 
    $n$~-~number of vertices/items; 
    $e$~-~number of edges;
    $d$~-~number of dimensions used to represent the constraints;
    $\zip$~-~optimal integer solution;
    $\zlp$~-~linear relaxation;
    $\#v, \#a$~-~number of vertices and arcs in the final arc-flow graph; 
    $\tpp$~-~time spent building the graph; 
    $\tlp$~-~time spent in the linear relaxation of the root node; 
    $\tip$~-~time spent in the branch-and-bound procedure; 
    $\nbb$~-~number of nodes explored in the branch-and-bound procedure; 
    $\ttot$~-~total run time in seconds. }
\end{threeparttable}
\end{table}

\subsubsection{Timetabling}
\label{sec:timetabling}

The timetabling problem has several variants and applications in many areas
(see, e.g., \citealt{TimetablesSA} and \citealt{Smith_hopfieldneural}).  In this
section, we consider the class-teacher-venue problem in which one has to find a
conflict-free timetable.  Suppose there are $c$ classes, $t$ teachers and $v$
venues.  Given the number of times each pair class-teacher must meet at each
venue, we want to find a timetable with zero clashes.  Classes, teachers and
venues cannot be chosen twice for the same time period.  These constraints are
represented by dimensions $\alpha$, $\gamma$ and $\delta$.  For each class $k$,
there is a dimension $\alpha^k$ of capacity 1.  For each teacher $k$, there is a
dimension $\gamma^k$ of capacity 1.  For each venue $k$, there is a dimension
$\delta^k$ of capacity 1.  Each requirement is a triplet (class $c_i$, teacher
$t_i$, venue $v_i$) that we represent by an item $i$ with weights
$\alpha^{c_i}_i = \gamma^{t_i}_i = \delta^{v_i}_i = 1$; 
the demand of each item is the demand of the corresponding requirement.  
In addiction,
there is a limit on the number of available time slots (that are represented by bins).  
This constraint is introduced in the objective function by searching for a solution 
that minimizes the number of periods.  This reduction can be seen as a graph coloring 
reduction to vector packing through clique constraints.  The set $S$ of valid patterns for
each time period is defined as follows:

\begin{equation}\label{eq:timetabling}
  A = \begin{bmatrix}
    \alpha_1^1 & \ldots & \alpha_m^1\\
    \svdots &        & \svdots\\
    \alpha_1^c & \ldots & \alpha_m^c\\
    \gamma_1^1 & \ldots & \gamma_m^1\\
    \svdots &        & \svdots\\
    \gamma_1^t & \ldots & \gamma_m^t\\    
    \delta_1^1 & \ldots & \delta_m^1\\
    \svdots &        & \svdots\\
    \delta_1^v & \ldots & \delta_m^v\\       
  \end{bmatrix}
  \qquad
  L = \begin{bmatrix}
    1\\
    \svdots \\    
    1\\
    1\\
    \svdots \\
    1\\
    1\\
    \svdots \\
    1\\        
  \end{bmatrix}
  \qquad 
  S = \{\XVEC \in \mathbb{N}_{0}^{m} : A\XVEC \leq L\}
\end{equation}

The arc-flow model was used to solve the ``hard timetabling'' instances from
\cite{ORLibrary}.  The hard classification comes from the fact that these
instances have been designed so that each class, teacher and venue is required
for each time period. Each instance has a solution with zero clashes that uses
no more than 30 periods (6 periods per day, 5 days per week).  All the instances
were solved quickly except hdtt8 that took a few hours to solve.
Table~\ref{tab:timetabling} shows the results for this problem.
The proposed method is very flexible and the 
solution of this problem shows that it is possible to model problems
beyond cutting and packing with very little effort.

\begin{table}[h!tbp]
\centering\begin{threeparttable}[b]
  \scriptsize 
  \caption{Results for timetabling.\label{tab:timetabling}}
  \begin{tabular}{lrrrrrrrrrrrr}
    \hline
    \up\down
    name & $t$ & $c$ & $v$ & $n$ & $m$ & $\#v$ & $\#a$ & $\tpp$ & $\tlp$ & $\tip$ & $\nbb$ & $\ttot$\\
    \hline
    \up    
    hdtt4 & 4 & 4 & 4 & 120 & 59 & 148 & 906 & 0.03 & 0.02 & 0.03 & 0 & 0.07\\
    hdtt5 & 5 & 5 & 5 & 150 & 88 & 644 & 4,221 & 0.14 & 0.15 & 0.29 & 0 & 0.58\\
    hdtt6 & 6 & 6 & 6 & 180 & 125 & 2,719 & 19,566 & 0.95 & 2.62 & 12.51 & 0 & 16.08\\
    hdtt7 & 7 & 7 & 7 & 210 & 154 & 11,140 & 82,725 & 5.59 & 47.62 & 1,697.31 & 16 & 1,750.52\\
    \down
    hdtt8 & 8 & 8 & 8 & 240 & 197 & 43,397 & 368,072 & 32.46 & 2,329.86 & 53,983.17 & 15 & 56,345.49\\
    \hline
  \end{tabular}
  \caption*{\footnotesize $t, c, v$~-~number of teachers, classes and
    venues.  $n$~-~number of items/requirements; $m$~-~number of different
    items/requirements; $\#v, \#a$~-~number of vertices and arcs in the final
    arc-flow graph; $\tpp$~-~time spent building the graph; $\tlp$~-~time spent
    in the linear relaxation of the root node; $\tip$~-~time spent in the
    branch-and-bound procedure; 
    $\nbb$~-~number of nodes explored in the branch-and-bound procedure;     
    $\ttot$~-~total run time in seconds. }
\end{threeparttable}
\end{table}

\FloatBarrier
\subsection{Bin packing and cutting stock}

Standard bin packing and cutting stock are one-dimensional vector packing problems
whose set $S$ of valid patterns is defined as follows:
\begin{equation}
  A = \begin{bmatrix}
    w_1 & \ldots & w_m  \\
  \end{bmatrix}
  \qquad
  L = \begin{bmatrix}
    W
  \end{bmatrix}
  \qquad 
  S = \{\XVEC \in \mathbb{N}_{0}^{m} : A\XVEC \leq L\}
\end{equation}

We used the arc-flow formulation to solve a large variety of bin packing
and cutting stock test data sets.
The results are summarized in Table~\ref{tab:bppcsp}.
\cite{ORLibrary} provides a bin packing test data set (BPP FLK) that was proposed by
\cite{springerlink:10.1007/BF00226291}.  This data set has two classes of instances:
uniform instances (uNNN), where items have randomly generated weights, and the
harder triplets instances (tNNN), where each bin in the optimal solution is
completely filled with three items.  Each of these is further divided into
subclasses of varying sizes.
We generated a cutting stock data set (CSP FLK) from this bin packing test data set 
by multiplying the demand of each item by one million;
note that in the class u1000 there are one thousand million items 
in each instance.
\cite{scholl} provides three data sets that were generated for
\cite{Scholl1997627}.  The first is composed of randomly generated instances
whose expected number of items per bin is not larger than 3.
The second test data set is composed of instances 
whose expected average number of items per bin is 3, 5, 7 or 9. 
Finally, the third test data set is composed of 10 difficult instances
with a total number of items $n = 200$ and bins of capacity $W = 100,000$. 
\cite{umetani} provides two large test data sets for cutting stock problems.
The first data set (Cutgen) is composed by 1800 randomly generated instances of 18
classes. These instances were generated using the problem generator proposed
by~\cite{Gau1995572}. The second data set (Fiber) was taken from a real application in
a chemical fiber company in Japan.
\cite{ShoenfieldHard28} provides the Hard28 test data set.  This data set is
composed of instances selected from a huge testing.  Among these 28 instances, 5
are non-IRUP, so the integer programming solver had to use branch-and-cut to
reduce the gap and prove optimality. The remaining instances are IRUP, yet very
hard for heuristics.
\cite{ESICUP} provides two test data sets collected from
\cite{ITOR:ITOR377} and \cite{springerlink:10.1007/BF01539705} (SCH/WAE and
WAE/GAU).
Finally, a cutting stock test data set was obtained from 1D-bar
relaxations of the two-dimensional bin packing test data set of \cite{Lodi1999}.

\begin{table}[h!tbp]
  \centering\begin{threeparttable}[b]
    \scriptsize 
    \caption{Results for the standard BPP/CSP.\label{tab:bppcsp}}
    \begin{tabular}{llrrrrrrrrrr}
      \hline
      \up\down
      data set & type & \#inst. & $m^{\max}$ & $n^{\max}$ & $W^{\max}$ & $\#v$ & $\#a$ & $\%v$ & $\%a$ & $\nbb$ & $\ttot$\\
      \hline
      \up
        BPP FLK & BPP & 160 & 203 & 1,000 & 1,000 & 116.26 & 3,803.94 & 49\% & 44\% & 1.46 & 0.63\\
        CSP FLK & CSP & 160 & 203 & $10^9$ & 1,000 & 115.38 & 3,844.81 & 49\% & 44\% & 0.00 & 0.60\\
        Fiber & CSP & 39 & 20 & 1,121 & 9,080 & 253.38 & 1,631.80 & 73\% & 71\% & 0.00 & 0.29\\
        Cutgen & CSP & 1,800 & 40 & 4,000 & 1,000 & 339.85 & 4,204.77 & 58\% & 51\% & 0.06 & 1.98\\
        1D-bar & CSP & 500 & 100 & 7,478 & 300 & 87.25 & 2,111.17 & 10\% & 7\% & 0.00 & 0.42\\
        Scholl & BPP & 1,210 & 350 & 500 & 100,000 & 294.81 & 13,669.59 & 34\% & 37\% & 0.04 & 17.42\\
        Hard28 & BPP & 28 & 189 & 200 & 1,000 & 789.46 & 27,284.00 & 19\% & 26\% & 102.57 & 29.68\\
        SCH/WAE & BPP & 200 & 49 & 120 & 1,000 & 210.11 & 3,553.57 & 70\% & 70\% & 0.00 & 0.62\\
        \down
        WAE/GAU & BPP & 17 & 64 & 239 & 10,000 & 6,235.06 & 128,212.76 & 33\% & 37\% & 2.59 & 1,641.09\\
      \hline
    \end{tabular}
    \caption*{\footnotesize 
        \#inst.~-~number of instances; 
        $m^{\max}$~-~maximum number of different items; 
        $n^{\max}$~-~maximum number of items;
        $W^{\max}$~-~maximum bin capacity; 
        $\#v, \#a$~-~average number of vertices and arcs in the final arc-flow graph; 
        $\%v, \%a$~-~average percentage of vertices and arcs removed by the graph compression method.
        $\nbb$~-~average number of nodes explored in the branch-and-bound procedure; 
        $\ttot$~-~average run time in seconds.         
    }
  \end{threeparttable}
\end{table}

\FloatBarrier
\subsection{Cardinality constrained bin packing and cutting stock}

One of the BPP variants is the cardinality constrained bin packing in which, in
addition to the capacity constraint, the number of items per bin is also
limited.  One of the variants of the CSP is cutting stock with cutting knife
limitation in which there is a limit on the number of pieces that can be cut
from each roll due to a limit on the number of knives.  BPP and CSP with
cardinality constraints can be seen as special cases of the 2-dimensional vector
packing problem.  In the 2-dimensional vector packing problem, 
there is a difficult problem in
each dimension, whereas, for the cardinality constrained BPP and CSP, the
problem is very easy in one of the dimensions: we just need to count the number
of items. The set $S$ of valid packing patterns for these problems is defined as
follows:

\begin{equation}
  A = \begin{bmatrix}
    w_1 & \ldots & w_m  \\
    1 & \ldots & 1 \\
  \end{bmatrix}
  \qquad
  L = \begin{bmatrix}
    W\\
    C\\
  \end{bmatrix}
  \qquad 
  S = \{\XVEC \in \mathbb{N}_{0}^{m} : A\XVEC \leq L\}
\end{equation}

Cardinality constrained bin packing is strongly NP-hard for any cardinality
larger than 2 (see, e.g., \citealt{Epstein:2011:IRM:1924502.1924506}); for
cardinality 2, the cardinality constrained bin packing problem can be solved in
polynomial time as a maximum non-bipartite matching problem in a graph where
each item is represented by a node and every compatible pair of items is connect
by an edge.

We solved using the arc-flow formulation
every instance from the bin packing and cutting stock data sets
with cardinalities between 2 and the minimum cardinality limit that allowed the optimum objective
value to be the same as the optimum without cardinality constraints.
Table~\ref{tab:card} summarizes the results for each data set.
The average run time in the 14,568 instances was 13 seconds
and 10,592 of these instances where solved in less than one minute.

\begin{table}[h!tbp]
  \centering\begin{threeparttable}[b]
    \scriptsize 
    \caption{Results for the cardinality constrained BPP/CSP.\label{tab:card}}
    \begin{tabular}{llrrrrrrrr}
      \hline
      \up\down
      data set & type & \#inst. & $C^{\max}$ & $\#v$ & $\#a$ & $\%v$ & $\%a$ & $\nbb$ & $\ttot$\\
      \hline
        \up
        BPP FLK & BPP & 320 & 3 & 58.58 & 1,700.37 & 86\% & 82\% & 1.44 & 0.35\\
        CSP FLK & CSP & 320 & 3 & 57.47 & 1,712.48 & 86\% & 83\% & 0.02 & 0.27\\
        Fiber & CSP & 279 & 12 & 83.23 & 525.35 & 94\% & 90\% & 0.00 & 0.07\\
        Cutgen & CSP & 7,299 & 18 & 253.15 & 2,459.78 & 93\% & 89\% & 0.06 & 1.26\\
        1D-bar & CSP & 1,415 & 8 & 59.81 & 965.96 & 82\% & 79\% & 0.05 & 0.25\\
        Scholl & BPP & 3,748 & 10 & 222.21 & 6,384.86 & 86\% & 87\% & 0.03 & 9.99\\
        Hard28 & BPP & 56 & 3 & 100.73 & 1,991.96 & 94\% & 94\% & 25.12 & 0.82\\
        SCH/WAE & BPP & 1,000 & 6 & 74.90 & 924.78 & 89\% & 90\% & 0.00 & 0.20\\
        \down
        WAE/GAU & BPP & 131 & 18 & 6,833.89 & 103,768.11 & 91\% & 90\% & 0.73 & 999.58\\
      \hline
    \end{tabular}
    \caption*{\footnotesize 
        \#inst.~-~number of instances; 
        $C^{\max}$~-~maximum cardinality limit; 
        $\#v, \#a$~-~average number of vertices and arcs in the final arc-flow graph; 
        $\%v, \%a$~-~average percentage of vertices and arcs removed by the graph compression method.
        $\nbb$~-~average number of nodes explored in the branch-and-bound procedure; 
        $\ttot$~-~average run time in seconds.         
    }
  \end{threeparttable}
\end{table}

The arc-flow formulation proved to work very well in cardinality constrained instances, 
for all values of cardinality.  
Graph compression reduces substantially the graph sizes
and usually leads to graphs with size comparable to the size without cardinality
constraints.
In fact, there are instances in which cardinality constraints
help to reduce the final graph size, thus leading to easier models.
For some of the instances we knew, by construction, that there would be a solution with at
most three items in each bin, but we were not aware of any good method in the
literature for solving the cardinality constrained BPP/CSP in general.  The arc-flow model
allowed us to solve the cardinality constrained BPP/CSP as easily as the standard BPP/CSP. 
Note that, when both capacity and cardinality constraints are active,
heuristic methods based on assignments tend to perform poorly
due to the difficulty in finding good solutions.

\FloatBarrier
\subsection{Cutting stock with binary patterns}

Cutting stock with binary patterns (0-1 CSP)
is a CSP variant in which items of each type 
may be cut at most once in each roll.
In this problem, pieces are identified by their types and 
some types may have the same width.
This problem usually appears as bar 
and slice relaxations of orthogonal packing problems 
(see, e.g., \citealt{Scheithauer1999} and ~\citealt{OPPBounds}).
Cutting stock with binary patterns can be modeled as
a vector packing problem with $m+1$ dimensions.
The binary constraints are introduced by $m$ binary dimensions
and the set $S$ of valid packing patterns for this problem can be defined as
follows:
\begin{equation}
  A = \begin{bmatrix}
    w_1 & w_2 & \ldots & w_m  \\
    1 & 0 & \ldots & 0\\
    0 & 1 & \ldots & 0\\
    \vdots & \vdots & \ddots & \vdots\\
    0 & 0 & \ldots & 1\\
  \end{bmatrix}
  \qquad
  L = \begin{bmatrix}
    W\\
    1\\
    1 \\
    \vdots \\
    1
  \end{bmatrix}
  \qquad 
  S = \{\XVEC \in \mathbb{N}_{0}^{m} : A\XVEC \leq L\}
\end{equation}

The arc-flow formulation was used to solve all the CSP instances
with binary patterns.
Table~\ref{tab:01csp} summarizes the results for each data set.
The average run time in the 2,499 instances was 5 seconds
and 98\% of the instances were solved in less than one minute.

\begin{table}[h!tbp]
  \centering\begin{threeparttable}[b]
    \scriptsize 
    \caption{Results for the 0-1 CSP.\label{tab:01csp}}
    \begin{tabular}{lrrrrrrrrrrrr}
      \hline
      \up\down
      data set & \#inst. & $m^{\max}$ & $n^{\max}$ & $W^{\max}$ & $\#v$ & $\#a$ & $\tpp$ & $\tlp$ & $\tgg$ & $\tip$ & $\nbb$ & $\ttot$\\
      \hline
      \up
        CSP FLK & 160 & 203 & $10^9$ & 1,000 & 451.47 & 4,313.02 & 0.46 & 0.14 & 0.20 & 0.65 & 0.07 & 1.25\\
        Fiber & 39 & 20 & 1,121 & 9,080 & 104.33 & 357.54 & 0.02 & 0.01 & 0.01 & 0.01 & 0.00 & 0.03\\
        Cutgen & 1,800 & 40 & 4,000 & 1,000 & 1,001.79 & 3,320.18 & 0.23 & 1.22 & 0.02 & 3.71 & 0.00 & 5.16\\
        \down
        1D-bar & 500 & 100 & 7,478 & 300 & 889.21 & 3,095.53 & 0.58 & 0.36 & 0.08 & 2.44 & 0.12 & 3.39\\
      \hline
    \end{tabular}
    \caption*{\footnotesize       
        \#inst.~-~number of instances;
        $m^{\max}$~-~maximum  number of different items;
        $n^{\max}$~-~maximum number of items;
        $W^{\max}$~-~maximum bin capacity;
        $\#v, \#a$~-~average number of vertices and arcs in the final arc-flow graph;               
        $\tpp$~-~average time spent building the graph;
        $\tlp$~-~average time spent in the linear relaxation of the root node; 
        $\tgg$~-~average time required to compute the linear relaxation using column-generation;
	    $\tip$~-~average time spent in the branch-and-bound procedure; 
	    $\nbb$~-~average number of nodes explored in the branch-and-bound procedure; 
	    $\ttot$~-~average run time in seconds.
      }
  \end{threeparttable}
\end{table}

Since this problem usually arises as a relaxation of
orthogonal packing problems, the linear relaxation
is usually enough as the gap between the linear relaxation
and the optimal integer solution is usually very small.
Gilmore-Gomory's model with column-generation 
is usually used to obtain the lower bounds for this problem.
During the column-generation procedure, it is usually necessary
to solve many knapsack problems and
the arc-flow graphs can also be used to solve these problems.
After building the arc-flow graph for a given instance, the pattern
with smallest reduced cost can be computed by evaluating 
$f(u) = \max_{(u',v,i) \in A:u'=u} (f(v)+c_i), f(\vT) = 0$,
where $c_i$ is the shadow price of the constraint associated 
with the demand of items of weight $w_i$ 
(the loss arcs are treated as items with a shadow price $c_0=0$).
The solution can be found in linear time using dynamic programming
and this method is much faster than solving 
each multi-constraint knapsack problem using \texttt{Gurobi}.
The Gilmore-Gomory run times ($\tgg$) presented in Table~\ref{tab:01csp} 
were obtained using the arc-flow graphs for solving the knapsack problems.
The time required to solve the linear relaxation
of the arc-flow models using interior point methods
is also presented ($\tlp$).
The column-generation approach presents better average run times
due to the small number of columns that is usually required to
find the optimal linear relaxation and the little time required
to solve each knapsack problem using the arc-flow graph.
Nevertheless, the time required to solve the linear relaxation of 
the arc-flow model is also very small.

\FloatBarrier
\subsection{Bin packing with conflicts}

The bin packing packing with conflicts (BPPC) is one of the 
most import bin packing variants.
This problem consists of the combination of 
bin packing with graph coloring.
In addiction to the capacity constraints,
there are compatibility constraints.
This problem can be solved as a vector packing
problem with $c+1$ dimensions, where $c$ is the number
of dimensions used to model conflicts.
The set $S$ of valid packing patterns for this problem can be defined as
follows:
\begin{equation}
  A = \begin{bmatrix}
    w_1 & \ldots & w_n  \\
    \alpha_1^1 & \ldots & \alpha_n^1\\
    \svdots &        & \svdots\\
    \alpha_1^c & \ldots & \alpha_n^c\\
  \end{bmatrix}
  \qquad
  L = \begin{bmatrix}
    W\\
    \beta^1\\
    \svdots \\
    \beta^c
  \end{bmatrix}
  \qquad 
  S = \{\XVEC \in \mathbb{N}_{0}^{n} : A\XVEC \leq L\}
\end{equation}

The conflicts can be modeled using any of the graph coloring
reductions to vector packing presented in Section~\ref{sec:color}.
In our experiments, we used degree constraints for 
modeling constraints in this problem.

In order to test the arc-flow formulation in the BPPC,
we used the data set proposed by \cite{Muritiba:2010:ABP:1838674.1838679}.
This data set was created from the bin packing data set of 
\cite{springerlink:10.1007/BF00226291}, adding
conflict graphs with several densities.
Tables~\ref{tab:bppc1} and \ref{tab:bppc2} summarize the results
for each class and density, respectively.
\cite{Muritiba:2010:ABP:1838674.1838679} solved
some of theses instances using a branch-and-price algorithm
under a time limit of 10 hours. 
\cite{sadykov:inria-00539869} solved all the instances,
including the open instances, within an one-hour time limit
using a branch-and-price algorithm.
As opposed to our method, both branch-and-price
algorithms of \cite{Muritiba:2010:ABP:1838674.1838679} and 
\cite{sadykov:inria-00539869} were specifically designed to solve
bin packing with conflicts.
The average run time of our method in the 800 instances was 2 minutes
and 80\% of these instances were solved in less than 1 minute.
The instances of class u1000 were the most difficult
to solve. Note that instances of this class have 1000 items
and, in the literature, instances with 200 items are already
considered difficult even in the one-dimensional case.
Nevertheless, no instance took longer than 50 minutes to be solved exactly.
In this problem, due to the high number of dimensions,
a large part of the run time is spent building the arc-flow graph.

\begin{table}[h!tbp]
  \centering\begin{threeparttable}[b]
    \scriptsize 
    \caption{Results for the BPP with conflicts.\label{tab:bppc1}}
    \begin{tabular}{lrrrrrrrrrrrr}
      \hline
      \up\down
      class & \#inst. & $n$ & $d$ & $d^{\max}$ & $\#v$ & $\#a$ & $\tpp$ & $\tlp$ & $\tgg$ & $\tip$ & $\nbb$ & $\ttot$\\
      \hline
      \up
		u120 & 100 & 120 & 84.96 & 121 & 359.99 & 3,012.12 & 0.18 & 0.06 & 0.12 & 0.18 & 0.07 & 0.43\\
		u250 & 100 & 250 & 175.21 & 251 & 1,167.75 & 9,393.41 & 2.21 & 0.37 & 0.63 & 1.29 & 0.00 & 3.86\\
		u500 & 100 & 500 & 350.48 & 501 & 3,486.90 & 27,440.63 & 32.40 & 1.87 & 3.82 & 10.88 & 0.00 & 45.16\\
		u1000 & 100 & 1,000 & 702.21 & 1001 & 11,316.90 & 88,323.32 & 643.68 & 13.24 & 28.35 & 104.79 & 0.00 & 761.72\\
		t60 & 100 & 60 & 42.22 & 61 & 99.69 & 693.14 & 0.03 & 0.01 & 0.03 & 0.02 & 0.00 & 0.06\\
		t120 & 100 & 120 & 84.01 & 121 & 298.36 & 2,627.20 & 0.20 & 0.05 & 0.11 & 0.27 & 3.70 & 0.52\\
		t249 & 100 & 249 & 174.42 & 250 & 1,045.20 & 10,300.08 & 2.82 & 0.32 & 0.50 & 1.65 & 3.12 & 4.78\\
		\down
		t501 & 100 & 501 & 352.26 & 502 & 3,796.35 & 36,809.01 & 53.90 & 2.56 & 2.64 & 17.50 & 0.00 & 73.95\\      
      \hline
    \end{tabular}
  \end{threeparttable}
\end{table}
\begin{table}[h!tbp]
  \centering\begin{threeparttable}[b]
    \scriptsize     
    \caption{Results for the BPP with conflicts (grouped by density).\label{tab:bppc2}}    
    \setlength{\tabcolsep}{7.2pt}
    \begin{tabular}{rrrrrrrrrrrr}    
      \hline
      \up\down
      density & \#inst. & $d$ & $d^{\max}$ & $\#v$ & $\#a$ & $\tpp$ & $\tlp$ & $\tgg$ & $\tip$ & $\nbb$ & $\ttot$\\
      \hline
      \up
		0\% & 80 & 1.00 & 1 & 113.47 & 12,288.20 & 0.28 & 0.24 & 4.79 & 1.24 & 0.00 & 1.76\\
		10\% & 80 & 69.76 & 222 & 814.02 & 20,034.35 & 18.35 & 0.62 & 4.77 & 4.01 & 0.00 & 22.98\\
		20\% & 80 & 140.60 & 420 & 2,077.11 & 24,327.59 & 75.96 & 1.38 & 7.66 & 11.45 & 0.00 & 88.79\\
		30\% & 80 & 211.31 & 632 & 3,540.24 & 29,554.71 & 161.19 & 2.77 & 11.05 & 50.21 & 8.61 & 214.18\\
		40\% & 80 & 280.51 & 818 & 4,951.05 & 35,266.04 & 233.93 & 4.73 & 9.67 & 70.75 & 0.00 & 309.40\\
		50\% & 80 & 350.02 & 1,001 & 5,862.62 & 39,442.55 & 227.11 & 6.55 & 3.89 & 17.02 & 0.00 & 250.68\\
		60\% & 80 & 351.00 & 1,001 & 4,456.75 & 29,639.79 & 126.87 & 4.03 & 2.04 & 10.46 & 0.00 & 141.36\\
		70\% & 80 & 351.00 & 1,001 & 3,020.04 & 19,390.99 & 56.31 & 1.95 & 0.93 & 4.06 & 0.00 & 62.32\\
		80\% & 80 & 351.00 & 1,001 & 1,643.91 & 10,220.38 & 16.17 & 0.73 & 0.35 & 1.37 & 0.00 & 18.28\\
        \down
		90\% & 80 & 351.00 & 1,001 & 484.70 & 3,084.05 & 3.10 & 0.09 & 0.09 & 0.16 & 0.00 & 3.35\\      		
      \hline
    \end{tabular}
    \caption*{\footnotesize 
        \#inst.~-~number of instances;
        $n$~-~number of items;
        $d$~-~average number of dimensions;
        $\#v, \#a$~-~average number of vertices and arcs in the final arc-flow graph;               
        $\tpp$~-~average time spent building the graph;
        $\tlp$~-~average time spent in the linear relaxation of the root node; 
        $\tgg$~-~average time required to compute the linear relaxation using column-generation;
	    $\tip$~-~average time spent in the branch-and-bound procedure; 
	    $\nbb$~-~average number of nodes explored in the branch-and-bound procedure; 
	    $\ttot$~-~average run time in seconds.     
      }
  \end{threeparttable}
\end{table}

\FloatBarrier
\subsection{Cutting stock with binary patterns and forbidden pairs}

Cutting stock with binary patterns and forbidden pairs (0-1 CSPC)
is a variant of cutting stock with binary patterns
that also includes compatibility constraints.
This problem usually appears as a relaxation of 
orthogonal packing problems (see, e.g., \citealt{OPPBounds}).
It can be modeled as a vector packing problem
with $c+m+1$ dimensions, where $c$ is the number
of dimensions used to model the conflicts.
The set $S$ of valid packing patterns for this problem can be defined as
follows:
\begin{equation}
  A = \begin{bmatrix}
    w_1 & w_2 & \ldots & w_m  \\
    \alpha_1^1 & \alpha_2^1 & \ldots & \alpha_m^1\\
    \svdots & \svdots &       & \svdots\\
    \alpha_1^c & \alpha_2^c & \ldots & \alpha_m^c\\
    1 & 0 & \ldots & 0\\
    0 & 1 & \ldots & 0\\
    \vdots & \vdots & \ddots & \vdots\\
    0 & 0 & \ldots & 1\\
  \end{bmatrix}
  \qquad
  L = \begin{bmatrix}
    W\\
    \beta^1\\
    \svdots\\
    \beta^c\\
    1\\
    1\\
    \vdots\\
    1
  \end{bmatrix}
  \qquad 
  S = \{\XVEC \in \mathbb{N}_{0}^{m} : A\XVEC \leq L\}
\end{equation}

In this problem, we also used degree constraints (see Section~\ref{sec:color}) 
to model the conflicts.
Since these constraints already guarantee that the corresponding
item cannot occur more than once in the same pattern,
the binary constraints for items with conflicts are discarded.
In order to test the behavior of the arc-flow model in this problem,
we created a 0-1 CSPC data set from the BPPC data set of 
\cite{Muritiba:2010:ABP:1838674.1838679}.
We assigned random values of demand between 1 and 100 to each item.
Tables~\ref{tab:bppc_cs1} and \ref{tab:bppc_cs2} summarize the
results for each class and density, respectively.
The average run time in the 800 instances was 6 minutes
and 72\% of these instances were solved in less than 1 minute.
Some instances took almost 3 hours to be solved and
1 instance took almost 20 hours.
Note that the graphs sizes for this type of problems
can be really large due to the high number of different items and dimensions;
in the class u1000, every instance has thousands of items of 1000 different
types and 1001 dimensions.

\begin{table}[h!tbp]
  \centering\begin{threeparttable}[b]
    \scriptsize 
    \caption{Results for the CSP with binary patterns and conflicts.\label{tab:bppc_cs1}}
    \setlength{\tabcolsep}{4pt}
    \begin{tabular}{lrrrrrrrrrrrrr}
      \hline
      \up\down
      class & \#inst. & $m$ & $n^{\max}$ & $d$ & $d^{\max}$ & $\#v$ & $\#a$ & $\tpp$ & $\tlp$ & $\tgg$ & $\tip$ & $\nbb$ & $\ttot$\\
      \hline
      \up
		u120 & 100 & 120 & 6,989 & 121 & 121 & 724.72 & 3,763.07 & 0.37 & 0.12 & 0.14 & 0.33 & 0.00 & 0.82\\
		u250 & 100 & 250 & 13,672 & 251 & 251 & 2,266.74 & 12,211.75 & 4.64 & 0.70 & 0.87 & 2.43 & 0.00 & 7.76\\
		u500 & 100 & 500 & 26,568 & 501 & 501 & 6,409.75 & 37,970.82 & 72.01 & 4.04 & 6.26 & 21.92 & 0.00 & 97.98\\
		u1000 & 100 & 1,000 & 52,492 & 1,001 & 1,001 & 19,171.67 & 126,689.21 & 1,518.38 & 28.63 & 50.38 & 1,198.08 & 48.58 & 2,745.09\\
		t60 & 100 & 60 & 3,567 & 61 & 61 & 110.98 & 709.32 & 0.04 & 0.01 & 0.03 & 0.03 & 0.00 & 0.08\\
		t120 & 100 & 120 & 6,856 & 121 & 121 & 322.66 & 2,637.34 & 0.32 & 0.05 & 0.10 & 0.20 & 0.44 & 0.57\\
		t249 & 100 & 249 & 13,578 & 250 & 250 & 1,102.75 & 10,311.04 & 5.17 & 0.35 & 0.63 & 2.29 & 0.00 & 7.81\\
        \down		
		t501 & 100 & 501 & 27,933 & 502 & 502 & 3,873.08 & 36,974.09 & 113.25 & 2.66 & 4.25 & 60.12 & 14.84 & 176.03\\
      \hline
    \end{tabular}
  \end{threeparttable}
\end{table}
\begin{table}[h!tbp]
  \centering\begin{threeparttable}[b]
    \scriptsize     
    \caption{Results for the CSP with binary patterns and conflicts (grouped by density).\label{tab:bppc_cs2}}  
    \setlength{\tabcolsep}{6.7pt}
    \begin{tabular}{rrrrrrrrrrrrrr}
      \hline
      \up\down
      density & \#inst. & $d$ & $d^{\max}$ & $\#v$ & $\#a$ & $\tpp$ & $\tlp$ & $\tgg$ & $\tip$ & $\nbb$ & $\ttot$\\
      \hline
	  \up
		0\% & 80 & 351 & 1,001 & 3,006.74 & 18,074.72 & 152.58 & 1.29 & 14.44 & 12.93 & 0.11 & 166.80\\
        10\% & 80 & 351 & 1,001 & 4,855.32 & 33,330.04 & 413.73 & 5.15 & 15.39 & 42.34 & 0.55 & 461.22\\
        20\% & 80 & 351 & 1,001 & 5,910.48 & 42,229.65 & 438.98 & 8.10 & 16.03 & 99.10 & 2.50 & 546.18\\
        30\% & 80 & 351 & 1,001 & 6,627.90 & 47,384.56 & 389.95 & 10.85 & 16.12 & 281.94 & 18.26 & 682.74\\
        40\% & 80 & 351 & 1,001 & 6,613.45 & 46,262.95 & 317.95 & 9.24 & 10.11 & 1,071.72 & 58.40 & 1,398.90\\
        50\% & 80 & 351 & 1,001 & 5,858.65 & 39,466.18 & 226.37 & 5.31 & 3.32 & 53.56 & 0.00 & 285.24\\
        60\% & 80 & 351 & 1,001 & 4,456.75 & 29,639.79 & 127.10 & 3.39 & 1.66 & 34.20 & 0.00 & 164.68\\
        70\% & 80 & 351 & 1,001 & 3,020.04 & 19,390.99 & 56.67 & 1.65 & 0.83 & 9.60 & 0.00 & 67.92\\
        80\% & 80 & 351 & 1,001 & 1,643.91 & 10,220.38 & 16.26 & 0.64 & 0.34 & 1.20 & 0.00 & 18.10\\
        \down
        90\% & 80 & 351 & 1,001 & 484.70 & 3,084.05 & 3.12 & 0.09 & 0.09 & 0.17 & 0.00 & 3.38\\		 
      \hline
    \end{tabular}
    \caption*{\footnotesize 
        \#inst.~-~number of instances;
        $m$~-~number of different items;
        $n^{\max}$~-~maximum number of items;
        $d$~-~average number of dimensions;
        $\#v, \#a$~-~average number of vertices and arcs in the final arc-flow graph;               
        $\tlp$~-~average time spent in the linear relaxation of the root node; 
        $\tgg$~-~average time required to compute the linear relaxation using column-generation;
	    $\tip$~-~average time spent in the branch-and-bound procedure; 
	    $\nbb$~-~average number of nodes explored in the branch-and-bound procedure; 
	    $\ttot$~-~average run time in seconds.   
    }
  \end{threeparttable}
\end{table}

\FloatBarrier
We are not aware of any effective exact method
from the literature for solving this problem.
In practice, Gilmore-Gomory's model with column-generation
is usually used to obtain good lower bounds for this problem
when it appears as a relaxation of other problems.
In Tables~\ref{tab:bppc_cs1} and \ref{tab:bppc_cs2}, 
we compare the run time of computing the linear
relaxation using column-generation (solving the knapsack problems
using the arc-flow graph) and the run time of computing
the linear relaxation of the arc-flow model using interior point methods.
Without considering the time required for building
the arc-flow graphs, both methods are very fast for computing
the lower bounds.
The time required to construct the graphs may be improved
by taking advantage of problem-specific characteristics;
for instance, \cite{Brandao01ArcFlow} use a single dimension
to represent binary constraints instead of $m$ binary dimensions.

\FloatBarrier
\subsection{Comparison with assignment-based formulations}

Assignment-based models are usually ineffective in practice
due to their symmetry and weak lower bounds.
The integer programming solvers usually include very powerful
heuristics which usually allow finding good solutions
for assignment-based models.
However, the optimal solution may be very difficult
to find.
Instances for which heuristics work well
can be solved really quickly using
assignment-based models;
on the other hand, when the heuristics dot not work well,
it is common to have small instances that take
days to be solved using an assignment-based formulation.
For instance, using \texttt{Gurobi} to solve the 
datased BPP FLK with an assignment-based formulation,
we were only able to solve 8 out of 160 instances
of this data set under an one-minute time limit;
increasing the time limit to 10 minutes, we were only able 
to solve more 5 instances.
Using the arc-flow model, we were able to solve
all the instances quickly with an average run
time less than 1 second.

Table~\ref{tab:comp} presents a comparison between
the percentages of instances solved 
using the assignment-based model and the arc-flow model
within an one-minute time limit.
The arc-flow model allows solving entirely many of the data sets
under an one-minute time limit
and the number of solved instances increases substantially
with slightly higher time limits.
There are few instances that take more than one hour
to be solved exactly using the arc-flow graph (see Section~\ref{sec:runtime}).

\begin{table}[h!tbp]
  \centering\begin{threeparttable}[b]
    \scriptsize 
    \caption{Comparison with assignment-based formulations.\label{tab:comp}}
    \begin{tabular}{rrrrrrrrrr}
      \hline
      \up\down
      & & & \multicolumn{3}{c}{standard} & & \multicolumn{3}{c}{card. constrained} \\
      \cline{4-6} \cline{8-10}
      \up\down
      data set & type & & \#inst. & \%K & \%AFG &  & \#inst. & \%K & \%AFG \\
      \hline
      \up
        BPP FLK & BPP &  & 160 & 5\% & 100\% &  & 320 & 53\% & 100\%\\
        CSP FLK & CSP &  & 160 & 0\% & 100\% &  & 320 & 0\% & 100\%\\
        Fiber & CSP &  & 39 & 36\% & 100\% &  & 279 & 90\% & 100\%\\
        Cutgen & CSP &  & 1,800 & 27\% & 100\% &  & 7,299 & 76\% & 100\%\\
        1D-bar & CSP &  & 500 & 23\% & 100\% &  & 1,415 & 49\% & 100\%\\
        Scholl & BPP &  & 1,210 & 41\% & 93\% &  & 3,748 & 77\% & 97\%\\
        Hard28 & BPP &  & 28 & 0\% & 86\% &  & 56 & 50\% & 100\%\\
        SCH/WAE & BPP &  & 200 & 57\% & 100\% &  & 1,000 & 92\% & 100\%\\
        \down
        WAE/GAU & BPP &  & 17 & 6\% & 18\% &  & 131 & 88\% & 61\%\\
      \hline
    \end{tabular}
  \end{threeparttable}
\end{table}
\begin{table}[h!tbp]
  \centering\begin{threeparttable}[b]
    \scriptsize   
    \setlength{\tabcolsep}{15pt}
    \begin{tabular}{rrrrr}
      \hline
      \up\down
      type & data set & \#inst. & \%K & \%AFG\\
      \hline
      \up        
        2CBP & 2CBP & 400 & 39\% & 76\%\\
        20d-VBP & 20CBP & 40 & 15\% & 78\%\\
        Coloring & Color & 4 & 25\% & 100\%\\
        Timetabling & Hard & 5 & 20\% & 60\%\\
        0-1 CSP & Cutgen & 1800 & 52\% & 98\%\\
        0-1 CSP & CSP FLK & 160 & 0\% & 100\%\\
        0-1 CSP & Fiber & 39 & 100\% & 100\%\\
        0-1 CSP & 1D-bar & 500 & 24\% & 98\%\\
        BPPC & BPPC & 800 & 3\% & 80\%\\
        \down 
        0-1 CSPC & BPPC\_CS & 800 & 0\% & 72\%\\
      \hline
    \end{tabular}    
    \caption*{\footnotesize 
      \#inst.~-~number of instances;
      \%K~-~percentage of instances solved using the assignment-based formulation
      in under an one-minute time limit;
      \%AFG~-~percentage of instances solved using the arc-flow formulation
      in under an one-minute time limit.}
  \end{threeparttable}
\end{table}

\FloatBarrier
\subsection{The importance of graph compression}
\label{sec:importance}

\begin{figure}[h!tbp]
  \caption{Graph size reduction (vertices).\label{red_vert}}
  \centering \includegraphics[scale=1]{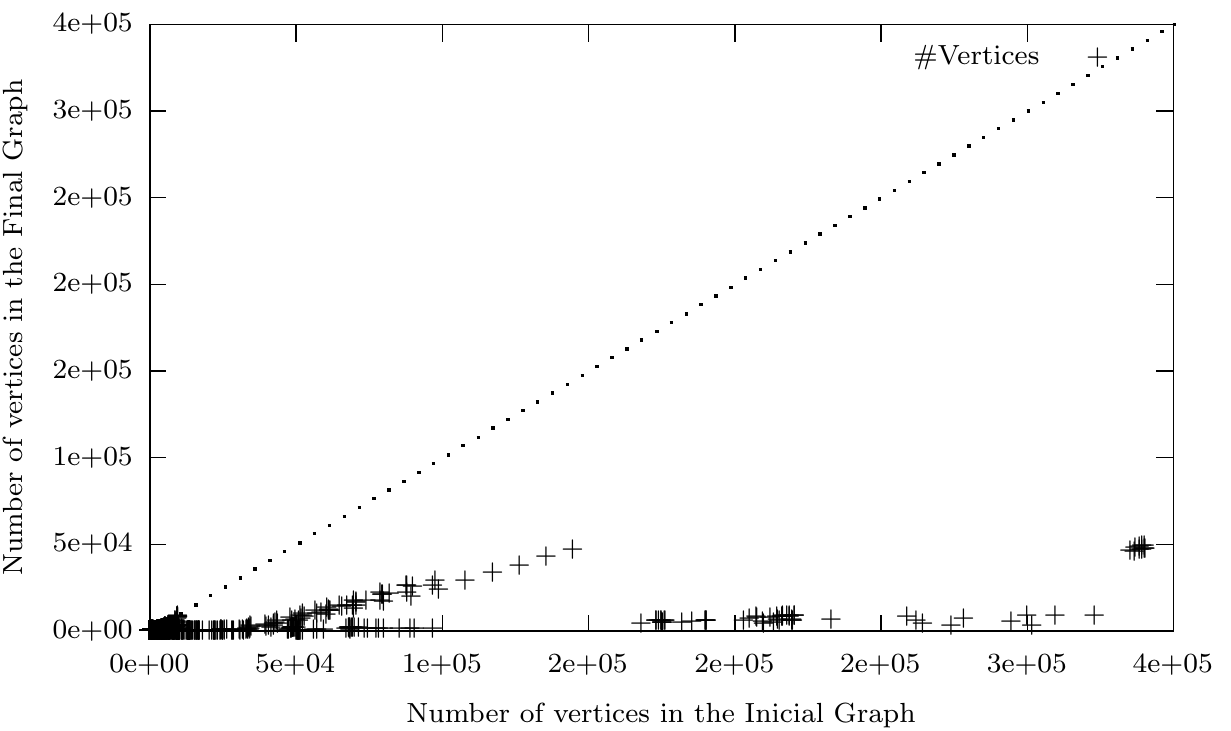}  
\end{figure}

\begin{figure}[h!tbp]
  \caption{Graph size reduction (arcs).\label{red_arcs}}
  \centering \includegraphics[scale=1]{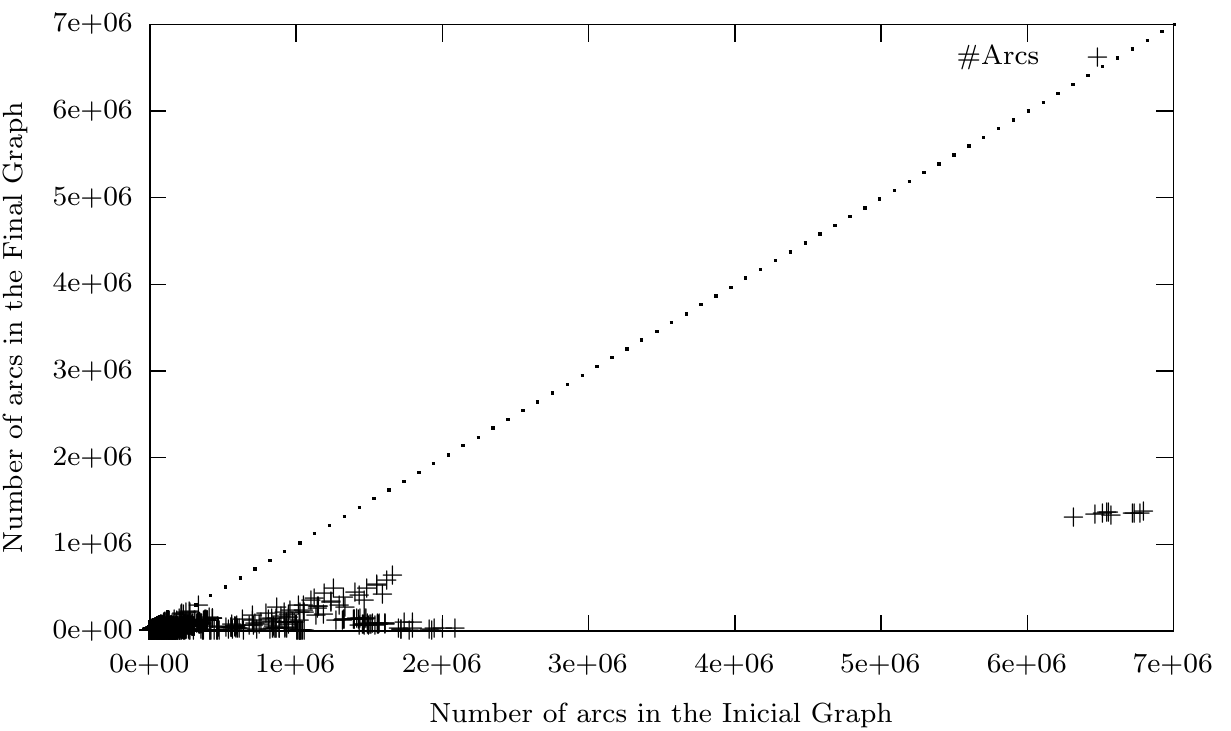}  
\end{figure}

Figures~\ref{red_vert} and \ref{red_arcs} show the relation
between the number of vertices and arcs before and after graph compression.  
We only consider instances with up to two dimensions (bin packing
and cutting stock with and without cardinality constraints, and
two-constraint bin packing) since the initial graph is 
usually huge for the remaining problems.
In these problems with up to two dimensions it is already possible
to see a remarkable graph size reduction.
Without graph compression it would be very difficult to
solve many of these instances within a reasonable amount of time.
For instance, arc-flow graphs for some two-constraint BPP instances
with 7 million arcs in the initial graph resulted in graphs
with approximately 1.5 million arcs. 
In the remaining problems, the compression ratios
are much higher. Problems whose initial graph
would be too large to fit in memory resulted
in reasonably small graphs due to graph compression;
it is common to obtain compression rates of hundreds of times.
Since we build the Step-3 graphs directly, we can
usually build the final arc-flow graph even when
the initial graph is too big to fit in memory.

\FloatBarrier
\subsection{Run time analysis}
\label{sec:runtime}

Using the proposed method, we solved sequentially 23,153 benchmark instances in
9 days, spending 33 seconds per instance, on average.  These benchmark
instances belong to several different problems.  
The same method was used to solve all the instances without 
any problem-specific adjustment.

\begin{figure}[h!tbp]
  \caption{Run time analysis (Gurobi).\label{runtimechart}}
  \centering \includegraphics[scale=1]{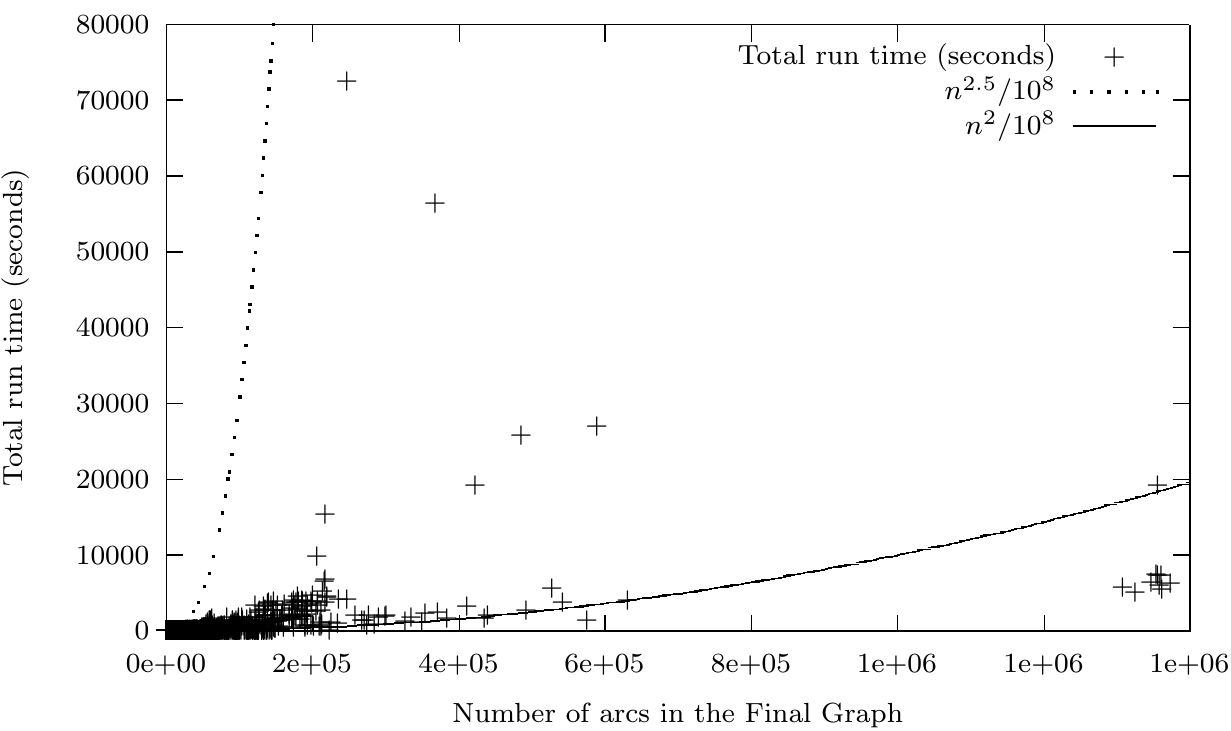}
\end{figure}

Figure~\ref{runtimechart} shows the relation between the number of arcs in
the final arc-flow graph and the total run time.  The two curves $n^2/10^8$ and
$n^{2.5}/10^8$ show an approximation of the run time (in seconds) of algorithms
with complexities $\Theta(n^2)$ and $\Theta(n^{2.5})$, with very low constant
factors. The large majority of the observed run times
appear between these two curves. Many of them are very close to the quadratic
run time, which is favorable since solving the arc-flow model is NP-hard.  The
few instances that lead to run times far from quadratic are mainly instances
where the number of items that fit in each bin is large (e.g., more than 10) and
hence the total number of patterns is huge.  Using heuristics, very good
solutions are usually found for these instances easily, since the waste tends to
be small.  Nevertheless, it may be very hard to find the optimal solution.
The arc-flow graphs for the 77 instances that were not
solved in a reasonable amount of time had millions
of arcs and the corresponding models are still too large
to be solved optimally using
current state-of-the-art mixed integer programming solvers
on a reasonable amount of time.

\begin{figure}[h!tbp]
  \caption{Run time analysis (GLPK).\label{runtimeglpkchart}}
  \centering \includegraphics[scale=1]{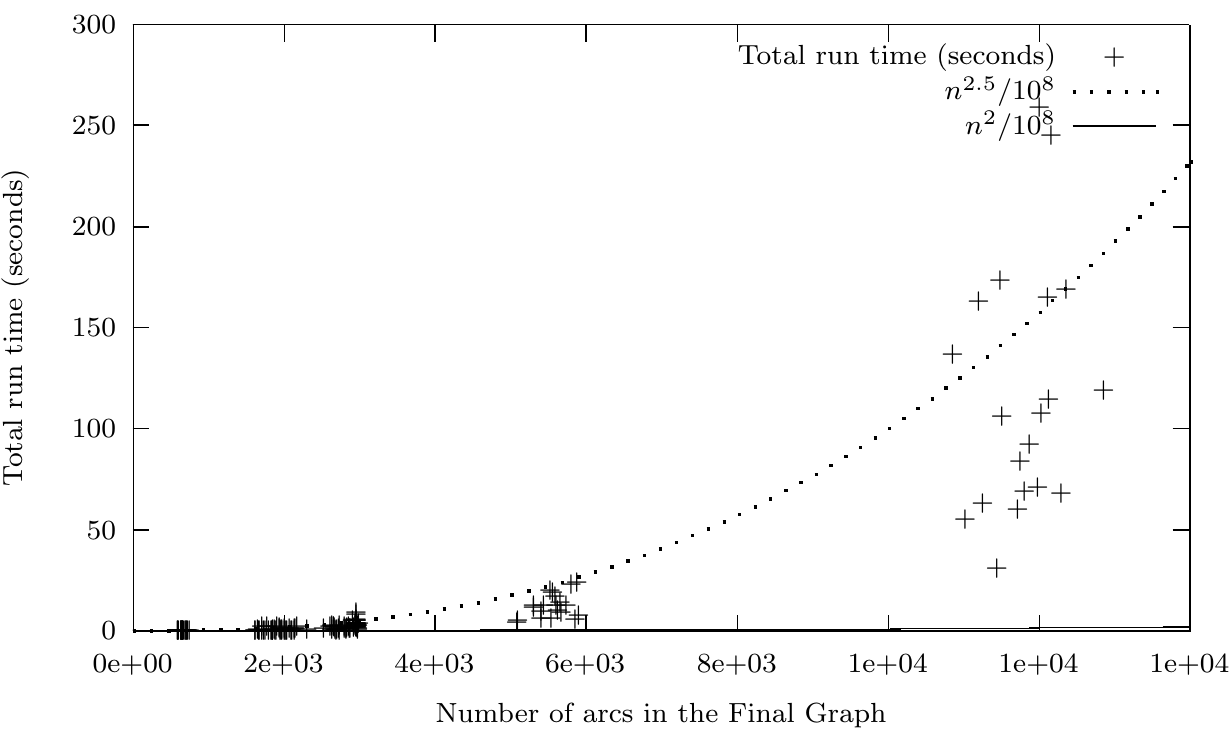}
\end{figure}

The bin packing data set BPP FLK of \cite{springerlink:10.1007/BF00226291} 
is one of the most used in the literature. In order to test
the effectiveness of our arc-flow model on non-commercial solvers,
this data set was solved using \texttt{GLPK}~4.43, an open-source mixed integer 
programming solver. 
Every instance was solved easily within a five-minute time limit.
Figure~\ref{runtimeglpkchart} shows the relation between the number of arcs in
the final arc-flow graph and the total run time.
Using \texttt{Gurobi} with an assignment-based formulation,
only 13 instances were solved within a ten-minute time limit,
showing that the formulation quality is extremely important even in the
one-dimensional case.

Very good results were also obtained using other non-commercial
MIP solvers such as \texttt{COIN-OR}.
Note that the non-commercial solvers are usually inferior to
commercial solvers and hence they may not be able to solve
the large models as easily as \texttt{Gurobi}.
Nevertheless, in practice, a non-commercial MIP solver
may be enough to solve instances that result in reasonably small models.
In the literature, most of the algorithms for solving
many of the problems solved here
are based on branch-and-price and usually rely
on very powerful MIP solvers such as IBM ILOG CPLEX. 
Our method is simple, effective and does not explicitly require
any particular MIP solver, though they may be necessary for solving large models.

\FloatBarrier
\section{Conclusions}
\label{sec:conclusions}

The method presented in this paper proved to be a very powerful tool
for solving several cutting and packing problems.  The model is
equivalent to Gilmore and Gomory's, thus providing a very strong
linear relaxation.  Nevertheless, it replaces column-generation by the
generation of a graph able to represent all the valid packing patterns
(one permutation of each pattern).
These are implicitly enumerated through the construction of a
compressed graph, which is proven to hold all the paths from the
source to the target that are required for determining the optimum
solution of the original problem.

This method can be used for solving several problems.
In this paper, we have dealt with
vector packing, graph coloring,
bin packing, cutting stock, cardinality constrained bin packing,
cutting stock with cutting knife limitation, cutting stock 
with binary patterns, bin packing with conflicts, 
and cutting stock with binary patterns and forbidden pairs.
The proposed general arc-flow formulation is equivalent to Gilmore-Gomory's
formulation and is very flexible.
The models produced by our method are very strong
and can be solved exactly by general-purpose mixed-integer programming
solvers; the model is very tight, and most of the solutions
are found in the root of the branch-and-bound tree.

The presented graph compression method is simple and proved to be very effective. 
In instances with many dimensions, it is common
to obtain graphs hundreds of times smaller than the initial graphs.
The combination of strong models with a reasonably small number
of variables and constraints makes this method very effective in practice.
Despite its simplicity and generality, the proposed method usually outperforms
more complex approaches such as branch-and-price algorithms.

Using the proposed method, we solved most of the known benchmark instances 
on a desktop computer, spending less than one minute per instance, on average.  
These benchmark instances belong to several different problems.  
The same method was used to solve all the instances without 
any problem-specific adjustment.
The linear relaxations are extremely strong in every problem
we considered.
The largest absolute gap we found in all the instances
from benchmark test data sets was 1.0027.

Several multi-constraint cutting and packing problems can be solved 
using the proposed method by means of 
reductions to vector-packing by defining a matrix of weights,
a vector of capacities and a vector of demands. 
Depending on the instance, it may or may not be possible to obtain
models of an acceptable size that can can be currently solved 
using any state-of-the-art mixed integer programming
solver; however,
our experiments show that it is very common to obtain 
reasonably small models on a large variety of problems,
thanks to the proposed graph compression technique.
We solved instances with up to 1,000 different items
and 1,000 constraints.
In the literature, even in the one-dimensional case, 
instances with 200 different items are already considered
very difficult.

\bibliographystyle{apalike} \bibliography{paper}
\end{document}